\documentclass[11pt]{article} % For LaTeX2e

\usepackage{fullpage}
\usepackage{hyperref}
\usepackage{url}
\usepackage{natbib}

\usepackage{booktabs}       % professional-quality tables
\usepackage{nicefrac}       % compact symbols for 1/2, etc.
\usepackage{xcolor}         % colors

\usepackage{amsmath, amsfonts,amsthm}
\usepackage{algorithm}
\usepackage{algorithmic}
\usepackage{mathtools}

\newcommand{\eqdef}{\coloneqq}

\newcommand{\sqn}[1]{\left\| #1 \right\|^2}
\newcommand{\sqnx}[1]{\left\| #1 \right\|^2_{\boldsymbol{\mathcal{X}}}}
\newcommand{\Exp}[1]{\mathbb{E}\!\left[ #1 \right]}
\newcommand{\Expc}[1]{\mathbb{E}\!\left[ #1\;|\;\mathcal{F}^t \right]}
\newcommand{\oma}{\omega_{\mathrm{av}}}

\usepackage{color}
\usepackage{xspace}
\usepackage{scalefnt}
\definecolor{mydarkred}{rgb}{0.75,0.0,0.0}
\newcommand{\algn}[1]{{\sf\scalefont{0.9}{#1}}\xspace}
\newcommand{\algno}{{\sf\scalefont{0.9}\dr{BiCo}\dg{LoR}}\xspace}
\newcommand\dr[1]{{\color{mydarkred}#1}}
\definecolor{mydarkgreen}{rgb}{0,0.5,0.1}
\newcommand\dg[1]{{\color{mydarkgreen}#1}}

\newcommand\al{\alpha}

\usepackage{colortbl}
\definecolor{mylgreen}{rgb}{0.87,1,0.82}
\definecolor{mygreen1}{rgb}{0,0.8,0}
\definecolor{gray}{rgb}{0.9,0.9,0.9}

\theoremstyle{plain}
\newtheorem{theorem}{Theorem}[section]

\newtheorem{corollary}[theorem]{Corollary}
\theoremstyle{definition}

\newtheorem{assumption}[theorem]{Assumption}
\theoremstyle{remark}

\title{\textbf{\dr{BiCo}\dg{LoR}: Communication-Efficient Optimization\\ with \dr{Bi}directional  \dr{Co}mpression and \dg{Lo}cal T\dg{r}aining}}

\author{Laurent Condat\qquad Artavazd Maranjyan\qquad  Peter Richtárik
\\
King Abdullah University of Science and Technology (KAUST)\\
Thuwal, %23955-6900, 
Kingdom of Saudi Arabia\\
\texttt{first.last@kaust.edu.sa} \\
}

\date{May 23, 2025}

\begin{document}

\maketitle

\begin{abstract}
 Slow and costly communication is often the main bottleneck in 
  distributed optimization, especially in federated learning 
  where it occurs over wireless networks.
 We introduce \algno, a communication-efficient optimization algorithm that combines two widely used and effective strategies: \dg{local training}, which increases computation between communication rounds,  and \dr{compression}, which encodes high-dimensional vectors into short bitstreams. While these mechanisms have been combined before,  compression
 has typically been applied only to uplink (client-to-server) communication, leaving the downlink (server-to-client) side unaddressed. In practice, however, both directions are costly.
  We propose
  \algno, the first algorithm to combine \dg{local training} with  \dr{bidirectional compression} using arbitrary unbiased compressors. 
  This joint design achieves accelerated complexity guarantees in both convex and strongly convex heterogeneous settings.
Empirically, \algno outperforms existing algorithms and establishes a new standard in communication efficiency.
\end{abstract}

{%\footnotesize
\renewcommand\baselinestretch{0}
\tableofcontents
\renewcommand\baselinestretch{1}
}

\section{Introduction}

Distributed computing has become pervasive across scientific disciplines. A prominent example is Federated Learning (FL), which enables collaborative training of machine learning models in a distributed manner \citep{kon16a,kon16,mcm17, bon17}. This rapidly growing field leverages data residing on remote devices, such as smartphones or hospital workstations. While FL must address challenges such as preserving privacy and resisting adversarial threats, \textbf{communication efficiency} is arguably the most critical issue \citep{kai19, li20,wan21}. 
Unlike centralized data-center setups, FL relies on parallel computation across client devices that must regularly exchange information with a distant central server. Communication typically occurs over bandwidth-limited and potentially unreliable networks, such as the internet or mobile networks. Consequently, communication remains the main bottleneck that hinders the widespread adoption of FL in large-scale consumer applications.

To mitigate communication overhead,
two main strategies have gained prominence: (1) \textbf{Local Training (\dg{LT})}, which reduces communication frequency by allowing multiple (stochastic) gradient steps to be performed locally before transmitting updates; and (2) \textbf{Communication Compression (\dr{CC})}, where updates are transmitted in a compressed form rather than as full-dimensional vectors. A review of the literature on LT and CC is presented in Section~\ref{secsota}.  

It is important to distinguish between \textbf{uplink communication (UpCom)} (clients to server) and \textbf{downlink communication (DownCom)} (server to clients). UpCom tends to be slower, as clients must upload distinct messages to the server that needs to decompress each of them, whereas in DownCom, all clients typically receive the same message simultaneously. Several factors can reinforce this asymmetry, such as limitations in cache size and aggregation capabilities at the server, as well as differences in network protocols or service provider configurations across internet or cellular systems. Many methods have been proposed that use compression for UpCom only, considering that DownCom is  cheap and can be neglected. This is not realistic, 
and in practical  settings both ways are
costly. For this reason, we focus on \textbf{Bidirectional Compression (\dr{BiCC})}, applied to both uplink and downlink communication.

We measure the \textbf{total communication (TotalCom)} cost, in number of bits, as a weighted sum of the UpCom and DownCom costs:
 \begin{equation}
 \mbox{TotalCom} = \mbox{UpCom} + \al.\mbox{DownCom},\label{eqtotcom}
  \end{equation}
  for some parameter $\al \geq 0$, typically in $(0,1]$. The case $\al=1$ corresponds to the symmetric regime in which UpCom and DownCom are equally costly, while  $\al=0$ corresponds to ignoring DownCom completely, which is not realistic. We focus in this work on settings where $\al$ is not negligible, so that downlink compression is required. The model \eqref{eqtotcom} has been proposed in \citet{con22cs,tamu23}, but counting the number of reals, not bits. 
  Since real numbers are represented with finite precision (typically 32-bit floats), and quantization  cannot reduce this size beyond a constant factor, see below, 
  both measures are proportional. Nevertheless, counting bits provides a more accurate measure.

\subsection{Unbiased compression}\label{secuc}

A common approach to reducing communication complexity in distributed learning is the use of lossy \dr{compression}. This involves applying a (potentially randomized) compression operator $\mathcal{C}: \mathbb{R}^d \to \mathbb{R}^d$ to a $d$-dimensional vector $x$, such that transmitting $\mathcal{C}(x)$, encoded as a short bit stream, is significantly more efficient than sending the full  vector $x$ (note that $\mathcal{C}(x)$ itself is the vector after decoding back to  $\mathbb{R}^d$).
Some compressors are unbiased; that is, $\Exp{\mathcal{C}(x)} = x$, where $\mathbb{E}[\cdot]$ denotes expectation. Others are biased \citep{bez20}. 
A widely used sparsifying compressor is \texttt{rand}-$k$, where $k \in [d] := \{1, \dots, d\}$. It randomly selects $k$ coordinates of $x$, scales them by $d/k$, and sets the remaining elements to zero. When the receiver knows which indices were selected—e.g., via a shared pseudo-random generator or by encoding them using a small overhead of at most $k\log_2 d$ bits—only those $k$ values need to be communicated, achieving a compression ratio of $d/k$. 
Aside from sparsification, which reduces the number of reals, quantization is another widely-used technique that reduces the number of bits needed to represent these reals. For instance, Natural compression rounds a real number to a power of 2 probabilistically and represents it using 9 bits, instead of 32 bits for  a full-precision float \citep{hor22}. 

To use compressors in iterative algorithms, we need to characterize them.  For this, we define, for every $\omega\geq 0$, the set $\mathbb{U}(\omega)$ 
of \textbf{unbiased compressors} $\mathcal{C}:\mathbb{R}^d\rightarrow \mathbb{R}^d$ with bounded
\textbf{relative variance} $\omega$; that is, 
such that
%\begin{equation}
$\Exp{\sqn{\mathcal{C}(x)-x}}\leq \omega \sqn{x}$, for every $x\in\mathbb{R}^d$. 
 Many practical compressors belong to this class  
 \citep{bez20,alb20,hor22,con22e}. 
Notably,  \texttt{rand}-$k\in \mathbb{U}\big(\frac{d}{k}-1\big)$, 
and Natural compression belongs to $\mathbb{U}(\frac{1}{8})$, so that its compression factor $32/9$ is almost free. Composing two compressors in  $\mathbb{U}(\omega_1)$ and $\mathbb{U}(\omega_2)$ yields a compressor in $\mathbb{U}(\omega)$ with $1+\omega=(1+\omega_1)(1+\omega_2)$ \citep{con22mu}. For instance \texttt{rand}-$1$ + Natural compression compresses a vector into $9 (+\log_2 d)$ bits, with $\omega = \frac{9d}{8}-1$. Importantly, there are lower bounds on the achievable compression level, and to compress vectors of $\mathbb{R}^d$ into $b$ bits, we have $\omega^{-1}\leq 4^{b/d}-1$, so that 
 $b(1+\omega)=\Omega(d)$ \citep{saf22,alb20}  \citep[Proposition 1]{yhe23}.
 
Moreover, given a collection $(\mathcal{C}_i)_{i=1}^n$ of  compression operators in $\mathbb{U}(\omega)$ for $\omega\geq 0$,  to characterize the relative variance after averaging their outputs, we introduce 
the constant $\oma\geq 0$ 
such that
\begin{equation}
\Exp{ \sqn{ \frac{1}{n}\sum_{i=1}^n \big(\mathcal{C}_i(x_i)-x_i\big)} } \leq \frac{\oma}{n} \sum_{i=1}^n \sqn{ x_i },
\label{eqbo}
\end{equation}
for every $(x_i)_{i=1}^n\in\big(\mathbb{R}^d\big)^n$. This is not an additional assumption, because \eqref{eqbo}  is satisfied with $\oma=\omega$. But $\oma$ can be much smaller than $\omega$. In particular, if the  $\mathcal{C}_i$ are mutually independent, \eqref{eqbo} is satisfied with $\oma=\frac{\omega}{n}$.

We introduce \algno, a novel randomized algorithm designed for communication-efficient  distributed optimization. It integrates \dg{LT} and \dr{BiCC} with arbitrary compressors  in $\mathbb{U}(\omega)$ for UpCom and $\mathbb{U}(\omega_s)$ for DownCom, for any $\omega,\omega_s\geq 0$. Variance reduction  with dual variables \citep{han19,gor202,gow20a} ensures convergence to the exact solution.

\subsection{Problem formulation} 

We consider the server-client model in which $n\geq 1 $ clients do computations in parallel and communicate in both directions with a server. 
We study distributed optimization problems of the form
\begin{equation}
\min_{x\in\mathbb{R}^d} \ \frac{1}{n}\sum_{i=1}^n f_i(x) +2f_s(x) + g(x),\label{eq1}
\end{equation}
where $d\geq 1$ is the model dimension, $f_i :\mathbb{R}^d\rightarrow \mathbb{R}$ is the private function of client $i\in[n]\eqdef\{1,\ldots,n\}$, $f_s$ is the private function of the server,  $g$ is a public function that is known by the server and  all clients. We suppose that 
 a solution $x^\star$ of \eqref{eq1} exists and we make the following assumptions.
\begin{itemize}
\item All functions $f_i$, $f_s$, $g$ are convex and $L$-smooth, for some $L> 0$ (a function $\phi$ is $L$-smooth if $\nabla \phi$ is $L$-Lipschitz continuous: for every $x,y$, $\|\nabla \phi(x)-\nabla \phi(y)\|\leq L \|x-y\|$. The norm is the Euclidean norm throughout the paper).  
\item  $f_i$, $f_s$, $g$ are  $\mu$-strongly convex, for some $\mu\geq 0$ ($\phi$ is $\mu$-strongly convex if $\phi-\frac{\mu}{2}\|\cdot\|^2$ is convex). $\mu$ can be zero in the general convex case. If $\mu>0$, the solution $x^\star$ of \eqref{eq1} exists and is unique, and we define the condition number $\kappa\eqdef \frac{L}{\mu}\geq 1$.
\end{itemize}

The problem \eqref{eq1} captures many important problems in machine learning, including   empirical risk minimization  \citep{sra11,shai_book}, and in many other fields. We aim to solve \eqref{eq1} efficiently in terms of communication, in the general \textbf{heterogeneous} setting; that is, the functions $f_i$, $f_s$, $g$  can be \emph{arbitrarily different}, without any similarity assumption. We prove accelerated communication complexity of \algno in both strongly convex and general convex cases. 

Thus, in the formulation \eqref{eq1}, the server is an additional $(n+1)$-th machine with its own loss function $f_s$, connected to the $n$ clients in a star-network. 
The idea of using an auxiliary dataset at the server representative of the global data distribution has been considered in several works, especially to correct for discrepancies induced by partial participation \citep{zhao18,yang21,yang23}. This is different from our setting, where $f_s$ can be very different from the $f_i$. 
The function $g$ is shared by all machines, so that Client $i$ makes calls to $\nabla f_i$ and $\nabla g$, and the server makes calls to $\nabla f_s$ and $\nabla g$. $g$ can also be the loss with respect to a small auxiliary dataset shared by all machines, but this is not necessarily the case. \emph{We stress that in the general convex case, $f_s$ and $g$ can be zero}.  In the strongly convex case, one can choose $f_s=g=\frac{\mu}{2}\sqn{\cdot}$. So, the template problem \eqref{eq1} is versatile and includes the minimization of $\sum_{i=1}^n f_i$ as a particular case. The reason why $g$ is introduced in \eqref{eq1} is that an estimate $y$ of a solution $x^\star$ is computed by all machines, and every client compresses the difference between its local model estimate and $y$. These differences tend to zero, which is key to obtaining a variance-reduced algorithm converging to a solution $x^\star$ exactly.

\subsection{Challenge and contribution}

This work addresses the following question: 
\emph{Can we combine \dg{LT} and \dr{BiCC} into a method that allows arbitrary unbiased compressors in UpCom and DownCom, provably benefits from the two techniques by exhibiting a state-of-the-art (SOTA) accelerated TotalCom complexity with nontrivial compression factors, and outperforms existing methods in practice?}\smallskip

Our new  algorithm \algno is the first to answer this complex question in the affirmative. 
This achievement overcomes significant challenges. For instance, with a standard application of random compression, the downlink and uplink compression errors  \emph{multiply} each other, 
as discussed in Section~\ref{seccompr}. In \algno, they are decoupled and only add up. 

\section{Related work}\label{secsota}

We review existing methods in the strongly convex case ($\mu>0$), as the study of the linear convergence rates provides important insights. We use the notation $\tilde{\mathcal{O}}(\cdot)=\mathcal{O}(\cdot \log \epsilon^{-1})$, where $\epsilon>0$ is the desired accuracy. We refer to \citet{tyu23} for a discussion of the general convex case.

\subsection{Local training}

Local Training (\dg{LT}) is a straightforward yet highly effective strategy. Instead of performing just a single (stochastic) Gradient Descent (GD) step between communication rounds, clients execute multiple steps locally. The core intuition is that these additional steps allow clients to transmit more informative updates, thereby reducing the total number of communication rounds required to achieve a target accuracy. \dg{LT} is a core component of the popular \algn{FedAvg} algorithm \citep{mcm17}, which is at the root of the immense success of FL. The idea of reducing the communication frequency was initially just heuristic, but a great deal of empirical evidence has demonstrated its practical effectiveness. \dg{LT} was analyzed first in the homogeneous data regime, or under restrictive assumptions such as bounded gradient
diversity \citep{had19}, then in the more realistic regime of heterogeneous data \citep{kha19,sti19, kha20a, li20a,woo20,gor21,gla22}. The more GD steps are made locally, the closer the local models get to the minimizers of the local functions $f_i$, which is not the desired behavior. This effect, called \emph{client drift}, 
has been quantified 
\citep{mal20}. The next class of methods, which includes \algn{Scaffold}  \citep{kar20},  \algn{S-Local-GD} \citep{gor21} and \algn{FedLin} \citep{mit21},  implemented variance reduction techniques to correct for client drift, so that a consensus is reached and every local model converges to the exact global solution. These methods are not accelerated, however.

More recently, a significant advancement was introduced by \citet{mis22} with \algn{Scaffnew}, the first  \dg{LT} algorithm that converges linearly 
with accelerated TotalCom complexity 
 $\tilde{\mathcal{O}}( d\sqrt{\kappa} )$. In \algn{Scaffnew}, communication occurs randomly 
 after each GD step with only a small probability 
$p$, resulting in an average of $1/p$  local steps between communication rounds. The optimal dependency $\sqrt{\kappa}$  \citep{sca19}
 is achieved when $p=1/\sqrt{\kappa}$.  \algn{Scaffnew} has then been extended  in  several ways \citep{mal22,maranjyan2023gradskip,con22rp,scafflix}.

\subsection{Compression}\label{seccompr}

A landmark development in the area of distributed algorithms using \dr{CC} is the variance-reduced algorithm  \algn{DIANA} proposed in 2019 \citep{mis19d}. It achieves linear convergence with uplink compressors in $\mathbb{U}(\omega)$ for any $\omega\geq 0$. Its iteration complexity, with communication at every iteration, is $\tilde{\mathcal{O}}\left(\left(1 + \frac{\omega}{n}\right)\kappa + \omega\right)$. 
So, with independent rand-1 compressors, its  UpCom complexity is $\tilde{\mathcal{O}} \left(\left(1 + \frac{d}{n}\right)\kappa + d\right)$, which significantly improves over  $\tilde{\mathcal{O}}(d\kappa )$ of standard GD when the number of clients $n$ is large. This has highlighted \dr{CC} as an acceleration mechanism.
 \algn{DIANA} has been extended in various directions, including support for stochastic gradients and partial participation \citep{hor22d,gor202,con22mu}. 
The algorithm  \algn{ADIANA} \citep{zli20}, based on Nesterov's accelerated GD,  has  iteration complexity
$\tilde{\mathcal{O}}\big(\big(1+\frac{\omega}{\sqrt{n}}\big)\sqrt{\kappa}+\omega\big)$ \citep{yhe23}.
For methods using independent uplink compressors  in $\mathbb{U}(\omega)$, the lower bound
$\tilde{\Omega}\big(\big(1+\frac{\omega}{\sqrt{n}}\big)\sqrt{\kappa}+\omega\big)$
on the  number of communication rounds has been established, so  \algn{ADIANA} is optimal in this sense.  This translates into the UpCom complexity $\tilde{\mathcal{O}}\big(\big(1+\frac{d}{\sqrt{n}}\big)\sqrt{\kappa}+d\big)$. 
Recently, linearly convergent algorithms using biased compressors have been proposed, such as \algn{EF21}
  \citep{ric21,fat21,con22e}, but the theoretical understanding of these methods remains less mature and their acceleration potential  is not clear. 
    
\dr{BiCC} is more complicated. A standard way to implement it is as follows: after UpCom, the server forms an updated model by aggregation of the received compressed vectors from the clients.  Then it compresses this model before DownCom to all clients for the next round. However, by doing so, the degradations due to uplink and then downlink compression pile up. A bidirectional extension of  \algn{DIANA} has been proposed, as part of the \algn{MURANA} framework \citep{con22mu}, with downlink compressors in $\mathbb{U}(\omega_s)$ for any $\omega_s\geq 0$. Its iteration complexity is $\tilde{\mathcal{O}}\!\left(\left(1 + \frac{\omega}{n}\right)\!(1+\omega_s)\kappa + \omega\right)$, as reported in Table~\ref{tab1}. We see that the uplink and downlink variances get multiplied with the $\left(1 + \frac{\omega}{n}\right)(1+\omega_s)$ dependence. Other algorithms, such as \algn{Artemis} \citep{phi20} and \algn{DORE} \citep{liu20} have been proposed, with the same complexity. \algn{MCM} \cite{phi21} has the slightly better complexity shown in Table~\ref{tab1}. 
Recently, \algn{EF21-P+DIANA} was proposed \citep{gru23}, extending \algn{DIANA} to \dr{BiCC} using error feedback at the server. Its complexity is $\tilde{\mathcal{O}}\left(\left(1 + \frac{\omega}{n}+\omega_s\right)\kappa + \omega\right)$, with better decoupled complexity depending on the sum $1 + \frac{\omega}{n}+\omega_s$ of the variances instead of their product. Lastly, \algn{2Direction} was introduced \citep{tyu23}, 
combining acceleration from momentum and decoupled \dr{BiCC}, achieving the SOTA complexity in number of rounds shown in Table~\ref{tab1}. With appropriate compressors, this gives a TotalCom complexity of $\tilde{\mathcal{O}}(d\sqrt{\kappa})$, which is the same, so neither worse nor better, than using no compression.  Also, we note that  \algn{2Direction} communicates full non-compressed vectors with a small probability. Ideally, a method would only communicate compressed vectors. \algno has this property and achieves the SOTA complexity $\tilde{\mathcal{O}}(d\sqrt{\kappa})$ as well. 
When $\alpha=0$ and no downlink compression is applied, 
neither \algn{2Direction}  nor \algno reverts to a known algorithm with unidirectional \dr{CC}, and their complexity is worse than \algn{ADIANA}, see the discussion in Section~\ref{secbounds}.

In a different area, \dr{BiCC} has been considered in the nonconvex Bayesian setting where compression consists of sampling from distributions \citep{egg25}.

\begin{table}[t]
\centering
\caption{Methods using arbitrary compressors $\mathcal{C}_i$ in $\mathbb{U}(\omega)$ for uplink and compressors $\mathcal{C}_s$ in $\mathbb{U}(\omega_s)$ for downlink communication, for arbitrary $\omega\geq 0$ and $\omega_s\geq 0$.  All compressors 
are independent.  
The $\tilde{\mathcal{O}}$ notation hides the log factors, in particular $\log(\epsilon^{-1})$, in  $\mathcal{O}(\cdot)$.\\
${}^{(a)}$ 2Direction requires communication of full non-compressed vectors with a small probability.\\
${}^{(b)}$ The reported complexity holds if $K$ satisfies the following conditions. For \algn{MURANA}: $K=\Omega\left(\frac{d}{n}\right)$. For \algn{MCM}: $K=\Theta(d)$ (no compression). 
For \algno: $K=\Omega\left(\frac{d}{\sqrt{\kappa}}\right)$.}
\label{tab1}
\begin{tabular}{l|l|l|l}
method&%\#\,
number of communication rounds
&TotalCom with rand-$K^{(b)}$\\
\hline
\algn{MURANA}&$\tilde{\mathcal{O}}\Big(\big(1+\frac{\omega}{n}\big)(1+\omega_s)\kappa+\omega\Big)$&
$\tilde{\mathcal{O}}\left( d\kappa\right)$\\
\algn{MCM}&$\tilde{\mathcal{O}}\left(\left(1+\frac{\omega}{n}+\omega_s^{3/2} +\frac{\sqrt{\omega}\omega_s}{\sqrt{n}}\right)\kappa\right)$&$\tilde{\mathcal{O}}\left(d\kappa\right)$\\
\algn{EF21-P+DIANA}&$\tilde{\mathcal{O}}\left((1+\frac{\omega}{n}+\omega_s)\kappa+\omega\right)$&$\tilde{\mathcal{O}}\left(d\kappa\right)$\\
\algn{2Direction}${}^{(a)}$&$
\tilde{\mathcal{O}}\left(\sqrt{(1+\omega)(1+\frac{\omega}{n}+\omega_s)\kappa}
+\omega+\omega_s\right)$&$\tilde{\mathcal{O}}\left(d\sqrt{\kappa}\right)$\\
\algno&$\tilde{\mathcal{O}}\left(\sqrt{(1+\omega+\omega_s)(1+\frac{\omega}{n}+\omega_s)\kappa}\right.$&$\tilde{\mathcal{O}}\left(d\sqrt{\kappa}\right)$\\
&\qquad$\left.+(1+\omega+\omega_s)(1+\frac{\omega}{n}+\omega_s)
\right)$
\end{tabular}
\end{table}

\subsection{Combining \dg{LT} and \dr{CC}}\label{secltcc}

It has proved difficult to combine \dg{LT} and \dr{CC} while keeping their benefits, namely acceleration from $\kappa$ to $\sqrt{\kappa}$ and an UpCom complexity with a better dependence on $d$ when $n$ is large. Early combinations, such as  \algn{Qsparse-local-SGD} \citep{bas20} and \algn{FedPAQ} \citep{rei20} fail to converge linearly. \algn{FedCOMGATE} \citep{had21} converges linearly  but in $\tilde{\mathcal{O}}(d \kappa)$. Random reshuffling, a technique that can be viewed as a kind of \dg{LT}, has been paired with \dg{CC} \citep{sad22, malinovsky2022federated}. The effective \dg{LT} mechanism of \algn{Scaffnew} has been combined with \dr{CC} in \algn{CompressedScaffnew}, achieving the UpCom complexity $\tilde{\mathcal{O}}\Big(\sqrt{d\kappa}+\frac{d\sqrt{\kappa}}{\sqrt{n}}+d \Big)$~\citep{con22cs}. It exhibits double acceleration with the $\sqrt{d}\sqrt{\kappa}$ dependence when $n$ is large. However, \algn{CompressedScaffnew} uses a specific linear compression technique based on random permutations of the coordinates. Recently, \algn{LoCoDL} was introduced \citep{con25loco}, successfully combining \dg{LT} in the spirit of \algn{Scaffnew} with \dr{CC} using arbitrary uplink compressors in $\mathbb{U}(\omega)$. 
Its UpCom complexity matches that of \algn{CompressedScaffnew}. \algn{ADIANA} has an even better complexity, that goes down to $\tilde{\mathcal{O}}(\sqrt{\kappa}+d) $ when $n$ is very large. Nevertheless, \algn{LoCoDL} consistently outperforms \algn{ADIANA}  in practice and can therefore be regarded as the SOTA in terms of UpCom efficiency. 

To the best of our knowledge, \algno is the first algorithm to combine the \dg{LT} mechanism of \algn{Scaffnew}, which yields $\sqrt{\kappa}$ acceleration, with \dr{BiCC} using arbitrary unbiased compressors. Similar to \algn{LoCoDL}, \algno  uses an additional function $g$ in the problem and a variable $y$ shared by all clients, with compression of the differences between the local variables $x_i$ and $y$. But it has notable differences: the server is an additional machine with its own function $f_s$, and during each communication round, the variables $x_i$ and $y$ are updated using information on $x_s$ received from the server, while $x_s$ itself is updated using information on the $x_i$ received from the clients. Crucially, these two updates are decorrelated, which is the key to obtain a decoupled TotalCom complexity, akin to \algn{EF21-P+DIANA} and \algn{2Direction}.

\begin{algorithm}[!t]
\small
	\caption{\algno}
	\label{alg2}
	\begin{algorithmic}[1]
		\STATE \textbf{input:}  stepsizes $\gamma,\eta,\eta_y,\rho,\rho_y>0$; 
		sequence of probabilities $(p_t)_{t\geq 1}$
		sparsification level $k\in[d]$;
		local initial estimates $x_1^0,\ldots,x_n^0,x_s^0,y^0 \in \mathbb{R}^d$, initial control variates 
		$u_1^0, \ldots, u_n^0,u_s^0,u_y^0 \in \mathbb{R}^d$ 
		such that $\frac{1}{n}\sum_{i=1}^n u_i^0 + 2u_s^0+u_y^0 =0$
		\FOR{$t=0,1,\ldots$}
				\FOR{$i=1,\ldots,n,s$, at clients and server in parallel,}
\STATE $\hat{x}_i^t\eqdef x_i^t -\gamma \nabla f_i(x_i^t) + \gamma u_i^t$
\STATE $\hat{y}^t\eqdef y^t -\gamma \nabla g(y^t) + \gamma u_y^t$ \ \ // the clients and server maintain identical copies of $y^t$, $u_y^t$
\ENDFOR
\STATE flip a coin $\theta^t\in\{0,1\}$, with $\mathrm{Prob}(\theta^t=1)=p_{t+1}$
\IF{$\theta^t=1$}
\STATE pick a subset $\Omega^t\subset[d]$ of size $k$ uniformly at random
\FOR{$i=1,\ldots,n$, at clients in parallel}
\STATE $c_i^t \eqdef \mathcal{C}_{i,\Omega^t}^t\big(\hat{x}_{i,\Omega^t}^t - \hat{y}_{\Omega^t}^t\big)$
\STATE send $c_i^t$ to the server
\STATE receive $c_s^t$ from the server
\STATE $x_{i,\Omega^t}^{t+1} \eqdef (1-\rho)\hat{x}_{i,\Omega^t}^t + \rho \big( c_s^t + \hat{y}^t_{\Omega^t}\big)$
\STATE $x_{i,[d]\backslash\Omega^t}^{t+1} \eqdef \hat{x}_{i,[d]\backslash\Omega^t}^t $
\STATE $y^{t+1}\eqdef \hat{y}^t+\rho_y c_s$
\STATE $u_i^{t+1} \eqdef u_i^t - \frac{p_{t+1}k \eta}{d\gamma}(c_i^t-c_s^t)$
\STATE $u_y^{t+1} \eqdef u_y^t + \frac{p_{t+1}k \eta_y}{d\gamma}c_s^t$
\ENDFOR
\STATE at server, in parallel  to steps 11--18:
\STATE\ \ \ \ $c_s^t \eqdef \mathcal{C}_{s,\Omega^t}^t\big(\hat{x}_{s,\Omega^t}^t - \hat{y}^t_{\Omega^t}\big)$
\STATE\ \ \ \ send $c_s^t$ to all clients
\STATE \ \ \ \ receive $(c_i^t)_{i=1}^n$ from the clients and aggregate $\bar{c}^t\eqdef \frac{1}{n}\sum_{i=1}^n c_i^t$
\STATE \ \ \ \ $x_{s,\Omega^t}^{t+1} \eqdef \left(1-\frac{\rho+\rho_y}{2}\right)\hat{x}_{s,\Omega^t}^t + \frac{\rho+\rho_y}{2}\hat{y}^t_{\Omega^t} + \frac{\rho}{2} \bar{c}^t $
\STATE \ \ \ \ $x_{s,[d]\backslash\Omega^t}^{t+1} \eqdef \hat{x}_{s,[d]\backslash\Omega^t}^t $
\STATE \ \ \ \ $y^{t+1}\eqdef \hat{y}^t+\rho_y c_s^t$
\STATE \ \ \ \ $u_s^{t+1} \eqdef u_s^t + \frac{p_{t+1}k \eta}{2d\gamma}\bar{c}^t -\frac{p_{t+1}k(\eta_y + \eta)}{2d\gamma}c_s^t$
\STATE \ \ \ \ $u_y^{t+1} \eqdef u_y^t + \frac{p_{t+1}k \eta_y}{d\gamma}c_s^t$
\ELSE
\STATE $x_i^{t+1}\eqdef \hat{x}_i^t\ \forall i\in[n]$,  $x_s^{t+1}\eqdef \hat{x}_s^t$, $y^{t+1} = \hat{y}^{t}$
\STATE $u_i^{t+1}\eqdef u_i^t\ \forall i\in[n]$,  $u_s^{t+1}\eqdef u_s^t$, $u_y^{t+1}\eqdef u_y^t$ 
\ENDIF
		\ENDFOR
	\end{algorithmic}
\end{algorithm}

\section{Proposed algorithm \algno}\label{secpro1}
 
The proposed stochastic primal--dual  method \algno  is shown as Algorithm~\ref{alg2}. At iteration $t$, Client $i$ computes  $\hat{x}_i^t$ by a GD step on its individual function $f_i$, corrected by a dual variable $u_i^t$ that learns $\nabla f_i(x^\star)$. It also computes $\hat{y}^t$ by a GD step on $g$. The server, as a $(n+1)$-th client, does the same using $f_s$ and $g$. Communication is triggered randomly with a small probability $p$. When it occurs, Client $i$ compresses  $\hat{x}_i^t-\hat{y}^t$ and sends this compressed difference to the server, which aggregates their average $\bar{c}^t$. Unlike in many algorithms, $\bar{c}^t$ is not sent back to the clients to update their local variables. Instead, it is only used by the server to update its local variables $x_s^t$ and $u_s^t$. This is the compressed difference $c_s^t=\mathcal{C}_s^t(\hat{x}_s^t-\hat{y}^t)$, sent by the server to all clients, which they use  to update their variables and $y$. So, UpCom and DownCom are independent and can be performed in parallel. We make the following assumption. 

\begin{assumption}[compressors in \algno]\label{ass2}
There exist $\omega,\omega_s\geq 0$ such that $\mathcal{C}_i^t\in\mathbb{U}(\omega)$ and $\mathcal{C}_s^t\in\mathbb{U}(\omega_s)$, for every $t\geq 0$, $i\in[n]$. The compressors $(\mathcal{C}_1^t,\ldots\mathcal{C}_n^t, \mathcal{C}_s^t)$ are independent from the $(\mathcal{C}_1^{t'},\ldots\mathcal{C}_n^{t'}, \mathcal{C}_s^{t'})$ if $t\neq t'$. Also, $\mathcal{C}_s^t$ is independent from the $(\mathcal{C}_i^t)_{i=1}^n$ for every $t\geq 0$. The $(\mathcal{C}_i^t)_{i=1}^n$ need not be mutually independent, this is characterized by $\oma$ in \eqref{eqbo}.
\end{assumption}

More technically, \algno works as follows.  We rewrite the problem \eqref{eq1} as
\begin{equation}
\min_{\mathbf{x}=(x_1,\ldots,x_n,x_s,y)} \ \frac{1}{n}\sum_{i=1}^n f_i(x_i) +2f_s(x_s) + g(y)\quad\mbox{s.t.}\quad\mathbf{D}\mathbf{x}=0,
\end{equation}
where $\mathbf{D}$ is a linear operator such that $\mathbf{D}\mathbf{x}=0$ if and only if $x_1=\cdots=x_n=x_s=y$. We refer to the Appendix for the vector notations and definitions. 
The key property of our design is that $\mathbf{D}$ is chosen in such a way that applying $\mathbf{D}$ and its adjoint $\mathbf{D}^*$ can be approximated using unbiased stochastic estimates, given that the clients receive the compressed vector $c_s$ from the server, and nothing else, and the server receives the compressed vectors $c_i$ from the clients.  Also, $(u_1^t,\ldots u_n^t,u_s^t,u_y^t)$ has to remain in the range of $\mathbf{D}^*$, which means that $\frac{1}{n}\sum_{i=1}^n u_i^t + 2u_s^t+u_y^t =0$, for every $t\geq 0$. This is why the idea that the server compresses the average of the $c_i^t$ and sends it back to the clients, instead of $c_s^t$, does not work, for instance. The operator norm of $\mathbf{D}$ is 2, implying that enabling \dr{BiCC} incurs a twofold slowdown of \algno without compression, relative to vanilla GD.
The iteration of \algno takes the form
\begin{equation*}
\left\lfloor\begin{array}{l}
\mathbf{\hat{x}}^t\eqdef \mathbf{x}^t-\gamma \nabla \mathbf{f}(\mathbf{x}^t)+\gamma \mathbf{D}^*\mathbf{u}^t\\
\mbox{flip a coin $\theta^t\in\{0,1\}$, with $\mathrm{Prob}(\theta^t=1)=p_{t+1}$}\\
\mbox{\textbf{if} }\theta^t=1:\ \ \mathbf{x}^{t+1}:\approx \mathbf{\hat{x}}^t-\boldsymbol{\rho}\mathbf{D}^*\mathbf{D}\mathbf{\hat{x}},\ \mathbf{u}^{t+1}:\approx \mathbf{u}^{t}- \frac{p\boldsymbol{\eta}}{\gamma}\mathbf{D}\mathbf{\hat{x}}\\
\mbox{\textbf{else}}:\ \ \mathbf{x}^{t+1}\eqdef \mathbf{\hat{x}}^t,\ \mathbf{u}^{t+1}\eqdef \mathbf{u}^t\\
%\mbox{\textbf{end if}}\\
\end{array}\right.,
\end{equation*}
where $:\approx $ means that an unbiased stochastic estimate of the right-hand side, built from the compressed vectors, is used for the update; see the Appendix  for precise definitions. The two stochastic estimates used for $\mathbf{x}$ and $\mathbf{u}$ are different. This is because an estimate of $\mathbf{D}^*\mathbf{D}\mathbf{\hat{x}}$ is needed to update $\mathbf{x}$, whereas an estimate of $\mathbf{D}\mathbf{\hat{x}}$ is needed to update $\mathbf{u}$, with $\mathbf{D}^*$ applied exactly to it. So, the estimate used for $\mathbf{x}$ is less noisy, with a variance that depends on $\oma+\omega_s$ instead of $\omega+\omega_s$.

In \algno, when communication happens, we allow the selection of a subset $\Omega^t$ of size $k\in[d]$ of coordinates to be processed. The other coordinates are updated as if $\theta^t=0$. For a vector $x\in\mathbb{R}^d$, its restriction $x_{\Omega}\in\mathbb{R}^d$ denotes $x$ with the coordinates not in $\Omega$ set to zero. Accordingly, $\mathcal{C}_\Omega(x_{\Omega})$ applies compression only to the subset of coordinates in $\Omega$ and returns a $k$-sparse vector with coordinates not in $\Omega$ set to zero. This approach allows to sparsify communication, but the same $k$ random coordinates are used for all vectors, in UpCom and DownCom. By contrast, we assume in Assumption~\ref{ass2} that if \texttt{rand}-$k$ compressors are used instead, achieving the same sparsification factor, the uplink and downlink compressors are independent. The interest of sparsifying via $k$ and not the compressors is that variance reduction of the compression error can be bypassed, leading to larger stepsizes $\eta$ and $\rho$.

\section{Convergence analysis and complexity of \algno}\label{seccc}

\subsection{Accelerated linear convergence in the strongly convex case}

\begin{theorem}[linear convergence of \algno]\label{theo1}
Suppose that $\mu>0$ and let $x^\star$ be the unique solution to \eqref{eq1}. 
In \algno, 
suppose that 
Assumption \ref{ass2} holds,  
$0<\gamma < \frac{2}{L}$, $p_t\equiv p\in (0,1]$ is constant, and
\begin{equation}
\rho=\rho_y=\frac{1}{2+\oma+2\omega_s}, \quad
\eta=\eta_y=\frac{1}{(1+2\omega+2\omega_s)(2+\oma+2\omega_s)}.
\end{equation}
For every $t\geq 0$, define the Lyapunov function
{\small\begin{align*}
\Psi^{t}&\eqdef  \frac{1}{\gamma}\!\left(\sum_{i=1}^n \sqn{x_i^t\!-\!x^\star} +2n \sqn{x_s^t\!-\!x^\star}+ n \sqn{y^t\!-\!x^\star}\right)+ \frac{d^2\gamma}{p^2k^2\eta}\!\left(\sum_{i=1}^n \sqn{u_i^t\!-\!u_i^\star} + n \sqn{u_y^t\!-\!u_y^\star}\right),
\end{align*}}%
where 
$u_y^\star \eqdef \nabla g(x^\star)$ and $u_i^\star \eqdef \nabla f_i(x^\star)$. Then 
\algno
converges linearly:  for every $t\geq 0$, 
\begin{equation}
\Exp{\Psi^{t}}\leq c^t \Psi^0,\quad\mbox{where}\quad
c\eqdef  \max\left((1-\gamma\mu)^2,(1-\gamma L)^2,1-\frac{p^2k^2\eta}{d^2}\right)<1.
\end{equation}
In addition,  for every $ i\in[n]$, $(x_i^t)_{t\in\mathbb{N}}$,  $(x_s^t)_{t\in\mathbb{N}}$ and $(y^t)_{t\in\mathbb{N}}$  converge  to $x^\star$, $(u_i^t)_{t\in\mathbb{N}}$ converges to $u_i^\star$, and $(u_y^t)_{t\in\mathbb{N}}$ converges to $u_y^\star$, 
almost surely.
\end{theorem}

Thus, \algno has the same rate 
$\max(1-\gamma\mu,\gamma L-1)^2$ as vanilla  GD, as long as $p^{-1}$ and the compression variances are below some threshold. 
The iteration complexity of \algno to reach $\epsilon$-accuracy, i.e.\ $\Exp{\Psi^{t}}\leq \epsilon$, 
with $\gamma=\Theta(\frac{1}{L})$, 
is
\begin{equation}
\mathcal{O}\left(\left(\kappa+\frac{d^2(1+\omega+\omega_s)(1+\oma+\omega_s)}{p^2k^2}\label{eq10}
\right)\log \frac{\Psi^0}{\epsilon}\right).
\end{equation}
With $k=d$ and $\omega_s=0$, this is the same complexity as \algn{LoCoDL}. 
The complexity in  number of communication rounds is $p$ times the iteration complexity, so 
the best value of $p$ balances the two terms in \eqref{eq10}; that is,
$
p\propto \min\left(\frac{d\sqrt{(1+\omega+\omega_s)(1+\oma+\omega_s)}}{k\sqrt{\kappa}},1\right)
$.
With this choice, the number of communication rounds is
\begin{equation}
\mathcal{O}\left(\left(\frac{d\sqrt{(1+\omega+\omega_s)(1+\oma+\omega_s)}}{k}\sqrt{\kappa}+\frac{d^2(1+\omega+\omega_s)(1+\oma+\omega_s)}{k^2}
\right)\log\frac{\Psi^0}{\epsilon}\right).\label{eqtcr}
\end{equation}
Assuming that the vectors compressed by $\mathcal{C}_i$ and $\mathcal{C}_s$ are encoded into $b$ and $b_s$ bits, respectively, the \textbf{TotalCom} complexity, as defined in \eqref{eqtotcom}, is $(b+\alpha b_s)$ times the complexity in \eqref{eqtcr}. As mentioned in Section~\ref{secuc},  $b(1+\omega)=\Omega(d)$, and the \texttt{rand}-$K$ compressor achieves this bound with $b=\Theta(K)$ (ignoring potential additional $K\log_2 d$ bits) and $1+\omega = \frac{d}{K}$. Thus, let us look at the TotalCom complexity with respect to $K\in[d]$ and $K_s\in[d]$, assuming that the compressors $\mathcal{C}_i$ and $\mathcal{C}_s$ are independent  \texttt{rand}-$K$ and   \texttt{rand}-$K_s$, respectively, and $k=d$. 
In addition, we assume that  $\alpha\leq 1$ and $K_s\geq K$. Then 
 TotalCom is
\begin{equation*}
\tilde{\mathcal{O}}\left((K+\alpha K_s)\left(
\sqrt{\frac{d}{K}\left(\frac{d}{nK}+\frac{d}{K_s}\right)}\sqrt{\kappa}+\frac{d}{K}\left(\frac{d}{nK}+\frac{d}{K_s}\right)\right)\right).
\end{equation*}
If in addition $nK\geq  K_s$, 
 this simplifies to 
$\tilde{\mathcal{O}}\left((K+\alpha K_s)\left(
\frac{d\sqrt{\kappa}}{\sqrt{KK_s}}
+\frac{d^2}{KK_s}\right)\right)$. 
In the unlikely case $\alpha\leq \frac{1}{d}$, with $K_s=d$ (no downlink compression since DownCom is very cheap) and $K=1$ (which implies $n\geq d)$, the complexity becomes $\tilde{\mathcal{O}}\left(\sqrt{d}\sqrt{\kappa}+d\right)$, which shows double acceleration with respect to $d$ and $\kappa$. In particular, when $\alpha=0$, this is the same complexity as \algn{CompressedScaffnew} 
and \algn{LoCoDL}  \citep{con25loco}. In the general case we are interested in, where $\alpha \in (0,1]$ is not tiny, notably the important case $\alpha=1$, we cannot expect acceleration with respect to $d$, see the discussion in Section~\ref{secbounds}. 
So, we suggest the following values.
\begin{corollary}\label{cor42}
In the conditions of Theorem~\ref{theo1}, suppose that $\alpha\in (0,1]$, $k=d$, $\gamma=\frac{1}{L}$, the $\mathcal{C}_i^t$ and $\mathcal{C}_s^t$ are independent \emph{\texttt{rand}-$K$} and  \emph{\texttt{rand}-$K_s$} compressors, respectively, with 
\begin{equation}
K_s = \left\lceil\frac{d}{\sqrt{\kappa}}\right\rceil,\quad K = \left\lceil\frac{\max\big(\alpha,\frac{1}{n}\big)d}{\sqrt{\kappa}}\right\rceil,\quad\mbox{and}\ \ p=\min\left(\frac{1}{\sqrt{\eta\kappa}},1\right).
\end{equation}
Then the analysis above applies, with $nK\geq K_s\geq K$ and $K\geq \alpha K_s$, and the TotalCom complexity of \algno is
{\small\begin{equation}
\mathcal{O}\left(\left(\frac{d\sqrt{K}\sqrt{\kappa}}{\sqrt{K_s}}
+\frac{d^2}{K_s}
\right)\log\frac{\Psi^0}{\epsilon}\right)=\mathcal{O}\left(d\sqrt{\kappa}\log\frac{\Psi^0}{\epsilon}\right).
\end{equation}}
\end{corollary}
An alternative is to make use of sparsification with the parameter $k$ in \algno:
\begin{corollary}\label{cor43}
In the conditions of Theorem~\ref{theo1}, suppose that $\alpha\in (0,1]$, $\gamma=\frac{1}{L}$, the compressors satisfy $\omega=\mathcal{O}(1)$, $\omega_s=\mathcal{O}(1)$ (e.g.\ quantization), 
$
k = \left\lceil\frac{d}{\sqrt{\kappa}}\right\rceil$, $ p=\min\left(\frac{d}{k\sqrt{\eta\kappa}},1\right)$. 
Then the TotalCom complexity of \algno is $\mathcal{O}\left(d\sqrt{\kappa}\log\frac{\Psi^0}{\epsilon}\right)$. 
\end{corollary}

We discuss the $\tilde{\mathcal{O}}(d\sqrt{\kappa})$ TotalCom complexity and explain why  it is unlikely to be improved in Section~\ref{secbounds}.

\subsection{Accelerated sublinear convergence in the general convex case}

For every $x\in\mathbb{R}^d$, $x'\in\mathbb{R}^d$, we define the  Bregman distance of a convex differentiable function $\phi$ as $\mathcal{D}_\phi(x,x')\eqdef \phi(x)-\phi(x') - \langle \nabla \phi(x'),x-x'\rangle \geq 0$. If $\phi$ is $L$-smooth, we have $\langle \nabla \phi(x)-\nabla \phi(x'),x-x'\rangle \geq \mathcal{D}_\phi(x,x') + \frac{1}{2L} \sqn{\nabla \phi(x)-\nabla \phi(x')}$.

\begin{theorem}[accelerated sublinear convergence of \algno]\label{theogc}
In \algno, %Algorithm \ref{alg2}, 
suppose that Assumption \ref{ass2} holds, $0<\gamma \leq \frac{1}{L}$, $k=d$,
\begin{equation}
\rho=\rho_y=\frac{1}{2+\oma+2\omega_s}, \quad
\eta=\eta_y=\frac{C}{(1+2\omega+2\omega_s)(2+\oma+2\omega_s)},
\end{equation}
for some constant $C\in(0,1)$, and that for every $t\geq 1$,
\begin{equation*}
p_t = \sqrt{\frac{b}{a+t}}\in (0,1],
\end{equation*}
for some  $b\geq \frac{1}{\eta}$ and $a\geq b-1$. Let $x^\star$ be a solution of \eqref{eq1}, and 
$u_y^\star \eqdef \nabla g(x^\star)$, $u_i^\star \eqdef \nabla f_i(x^\star)$ for every $i\in[n]$.
Then \algno converges sublinearly: 
for every $\epsilon>0$, by choosing $\tilde{t}$ uniformly at random in $\{0,\ldots,T-1\}$, where
{\small
\begin{align*}
T\!&\!\eqdef\left\lceil \!\frac{1}{2\epsilon}\!\left( \!\frac{1}{\gamma}\sum_{i=1}^n \sqn{x_i^0-x_i^\star}\! +\!  \frac{2n}{\gamma} \sqn{x_s^0-x_s^\star}\!+\!\frac{n}{\gamma} \sqn{y^0-y^\star}
\!+\!\frac{\gamma a}{\eta b}\sum_{i=1}^n\sqn{u_i^0-u_i^\star}\!+\!\frac{n\gamma a}{ \eta b}\sqn{u_y^0-u_y^\star}\!\right)\!\right\rceil\!,
\end{align*}}
we have\vspace{-1mm}
\begin{equation}
\Exp{\sum_{i=1}^n  \mathcal{D}_{f_i}\big(x_i^{\tilde{t}},x^\star\big)
+2n\mathcal{D}_{f_s}\big(x_s^{\tilde{t}},x^\star\big) +n \mathcal{D}_{g}\big(y^{\tilde{t}},x^\star\big)
}\leq \epsilon.\label{eqe2}
\end{equation}
and\vspace{-1mm}
\begin{align}
&\Exp{\frac{1}{\gamma}\sum_{i=1}^n 
\sqn{x_i^{\tilde{t}}-x_s^{\tilde{t}}}}
=\mathcal{O}\left(\sqrt{\epsilon}\right),\quad
\Exp{\frac{n}{\gamma}\sqn{x_s^{\tilde{t}}-y^{\tilde{t}}}}=\mathcal{O}\left(\sqrt{\epsilon}\right)
\end{align}
Moreover, the expectation of the number of communication rounds over the first $T\geq 1$ iterations is 
$
\sum_{t=1}^{T} p_t = \Theta(\sqrt{T})
$, 
so that \eqref{eqe2} is achieved with $\Theta(\sqrt{T})=\Theta\left(\frac{1}{\sqrt{\epsilon}}\right)$ communication rounds
(we refer to the proof in the Appendix for the constants in $\mathcal{O}$). Moreover, if $\gamma=\Theta(\frac{1}{L})$, $x^0\eqdef x_1^0=\cdots=x_n^0=x_s^0=y^0$ and $u_1^0=\nabla f_1(x^0),\ldots, u_n^0=\nabla f_n(x^0)$, $u_y^0=\nabla g(x^0)$, $\sqrt{T}=\Theta\left(
\sqrt{\frac{L}{\epsilon}}\|x^0-x^\star\|\right)$. This is the same accelerated complexity as \algn{2Direction} \citep{tyu23}.
\end{theorem}

\section{Discussion of the $\tilde{\mathcal{O}}(d\sqrt{\kappa})$ TotalCom complexity}\label{secbounds}

Let us comment on existing lower and upper bounds for the TotalCom complexity in the strongly convex setting.

As discussed in Section~\ref{seccompr}, 
 for methods using uplink compression only with independent compressors $\mathcal{C}_i^t$ in $\mathbb{U}(\omega)$, the lower bound, achieved by \algn{ADIANA}, for the  number of communication rounds is $\tilde{\Omega}\big(\big(1+\frac{\omega}{\sqrt{n}}\big)\sqrt{\kappa}+\omega\big)$.  This translates into an UpCom (or equivalently TotalCom with $\alpha=0$) complexity  of $\tilde{\mathcal{O}}\big(\big(1+\frac{d}{\sqrt{n}}\big)\sqrt{\kappa}+d\big)$. \algno, like \algn{2Direction}, achieves the same complexity $\tilde{\mathcal{O}}\big(\sqrt{d}\sqrt{\kappa}+d\big)$ as \algn{LoCoDL}  if $n\geq d$ and $\alpha=0$. This is worse than \algn{ADIANA}, but \algn{LoCoDL}  outperforms \algn{ADIANA}  in practice \citep{con25loco}. In any case, there remains a theoretical gap between unidirectional and bidirectional \dr{CC}.
 
Without downlink compression, or little compression with $\omega_s=\mathcal{O}(1)$ such as quantization solely, 
by choosing the $\mathcal{C}_i^t$ as independent \texttt{rand}-$K$ compressors for $K=\max( \frac{d}{n},1,d\alpha)$ (which means no uplink compression if $\alpha=1$), \algno achieves the same TotalCom complexity $\tilde{\mathcal{O}}\big(\frac{d\sqrt{\kappa}}{\sqrt{n}}+\sqrt{d\kappa}+d+\sqrt{\alpha}d\sqrt{\kappa}\big)$ as \algn{CompressedScaffnew} \citep{con22cs}. For not-so-small values of $\alpha$, this reverts to $\tilde{\mathcal{O}}(d\sqrt{\kappa})$. Moreover, this is obtained by essentially disabling compression and using \dg{LT} only. This is not in the spirit of what we want to achieve, which is the best TotalCom complexity with nontrivial levels of compression in both ways. In Corollaries \ref{cor42} and \ref{cor43}, \algno has complexity $\tilde{\mathcal{O}}(d\sqrt{\kappa})$ with large levels of compression, even when $\alpha=1$.
Still, there is a gap here too, between the regime $\alpha=0$ where acceleration with respect to $d$ is possible thanks to compression, and $\alpha=1$, where the TotalCom complexity is $\tilde{\mathcal{O}}(d\sqrt{\kappa})$ and compression is at best harmless. 
Nesterov's \algn{accelerated GD} and \algn{Scaffnew} \citep{mis22} have this complexity, and they don't use compression.

 It has been shown that without assuming independence of the $\mathcal{C}_i^t$, a lower bound for the number of communication rounds is $\tilde{\Omega}\big((1+\omega)\sqrt{\kappa}\big)$, which gives an UpCom complexity of $\tilde{\Omega}(d\sqrt{\kappa})$ \citep{he23l}. This complexity is achieved by Nesterov's \algn{accelerated GD} and \algn{Scaffnew} \citep{mis22}, which do not use compression. This may indicate that many independent compressors run in parallel are required to hope for a decrease of the dependence with respect to $d$. We may consider the idea that the server sends different messages compressed with  independent compressors, instead of using a single compressor $\mathcal{C}_s^t$. However, a negative result has been established, in the different nonconvex setting though: any method in which the server sends a compressed vector to each client, possibly obtained using different compressors in $\mathbb{U}(\omega_s)$, requires at least $\Omega\big( (1+\omega_s)L/\epsilon\big)$ rounds to find a stationary point 
 \citep[Theorem 3.1]{gru24}; see also \citet{hua22}. This suggests that there is little hope of improving the downlink dependence on $\omega_s$ as $n$ increases.

\begin{figure}[t]
    \centering
    \begin{tabular}{cc}
    \includegraphics[width=0.47\linewidth]{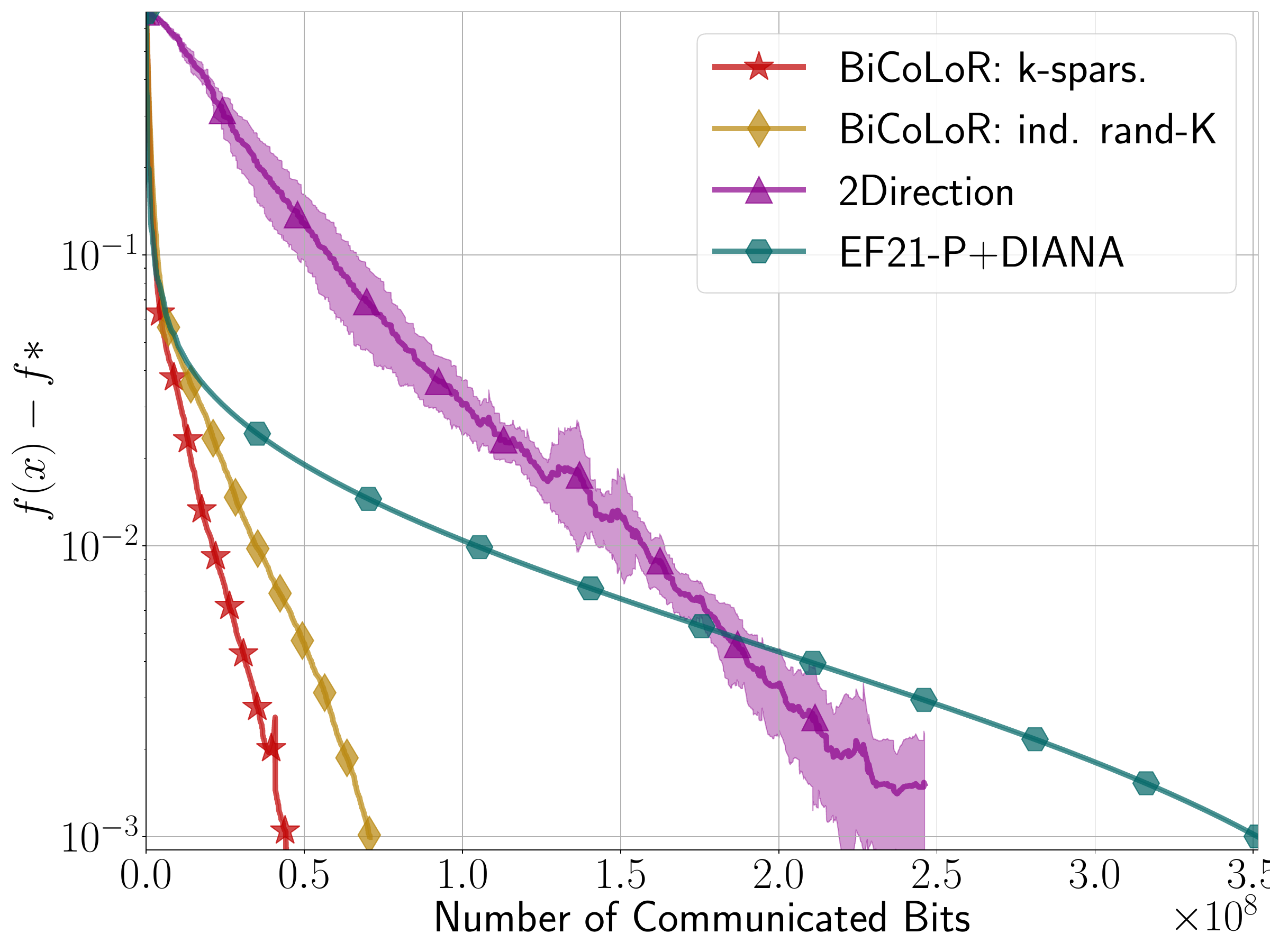}
    & 
    \includegraphics[width=0.47\linewidth]{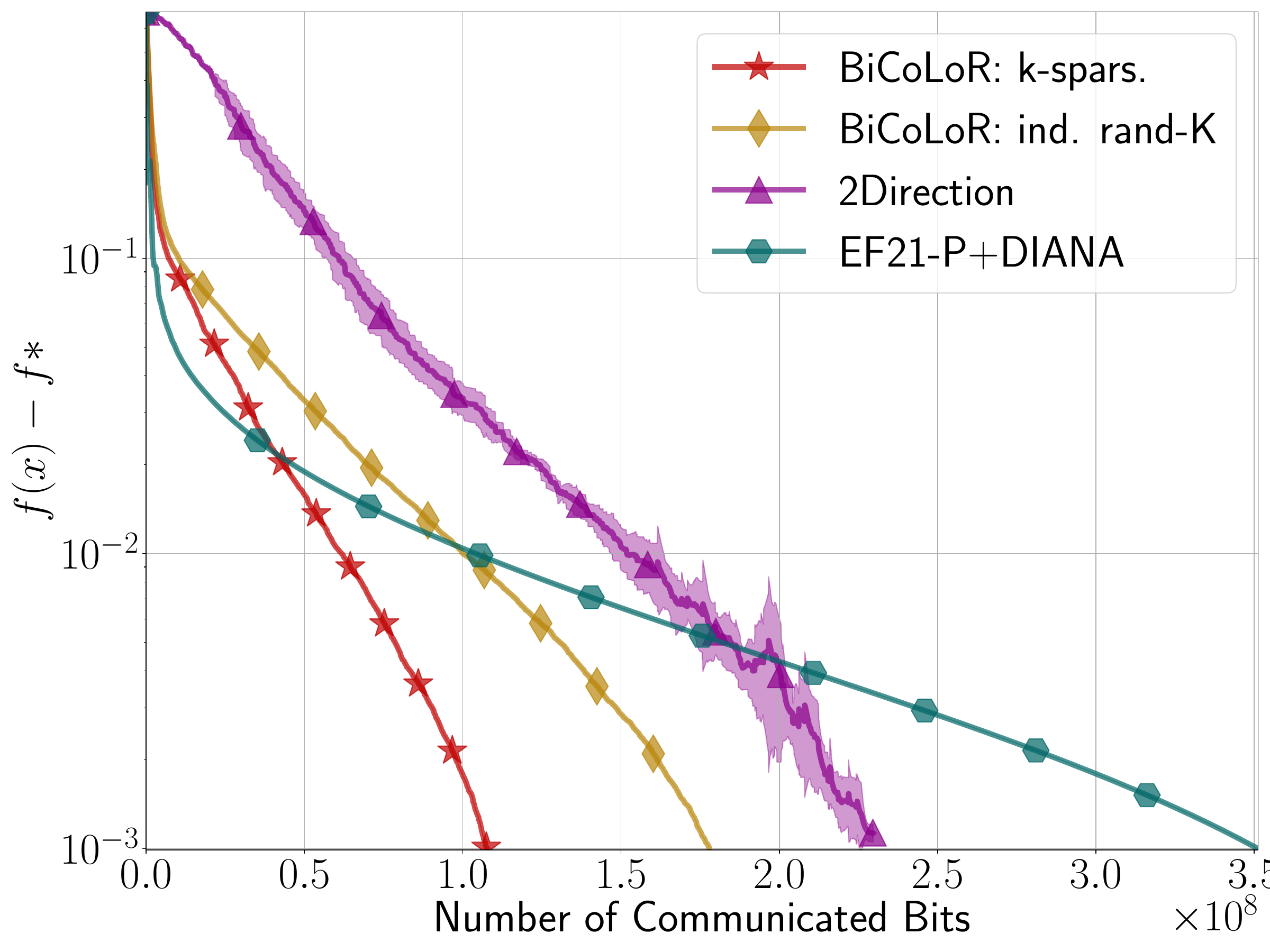} \\
    (a) $n=10$ &
    (b) $n=100$
    \end{tabular}
    \caption{
       Logistic regression on the \texttt{real-sim} dataset.
       The compression scheme combines sparsification with $K = 1000$ and Natural Compression.
    }
    \label{fig:real_sim}
\end{figure}

\section{Experiments}\label{secex}

\subsection{Strongly convex case}\label{secssc}

We evaluate our proposed method, \algno, against \algn{2Direction} and \algn{EF21-P+DIANA} on a regularized logistic regression problem of the form \eqref{eq1}.
The loss function of Client $i$ is
\begin{equation}\label{eqlogr}
    f_i(x)=\frac{1}{m} \sum\limits_{s=1}^{m} \log\!\Big(1\!+\!\exp \left(-b_{i,j} a_{i,j}^{\top} x\right)\!\Big)\!+\frac{\mu}{2}\|x\|^{2},
\end{equation}
and we take $f_s = g = \frac{\mu}{2}\|x\|^{2}$.
For the other algorithms, which do not use $f_s$ and $g$, we replace $\mu$ by $4\mu$ in the functions $f_i$, so that the problem solved by the different algorithms is exactly the same. 
In \eqref{eqlogr}, $n$ is the number of clients, $m$ is the number of data points per client, $a_{i,j} \in \mathbb{R}^{d}$ and $b_{i,j} \in \{-1,+1\}$ are data samples, and $\mu$ is set so that the condition number is $\kappa = 4.10^6$.

\begin{figure}
    \centering
    \begin{tabular}{cc}
    \includegraphics[width=0.47\linewidth]{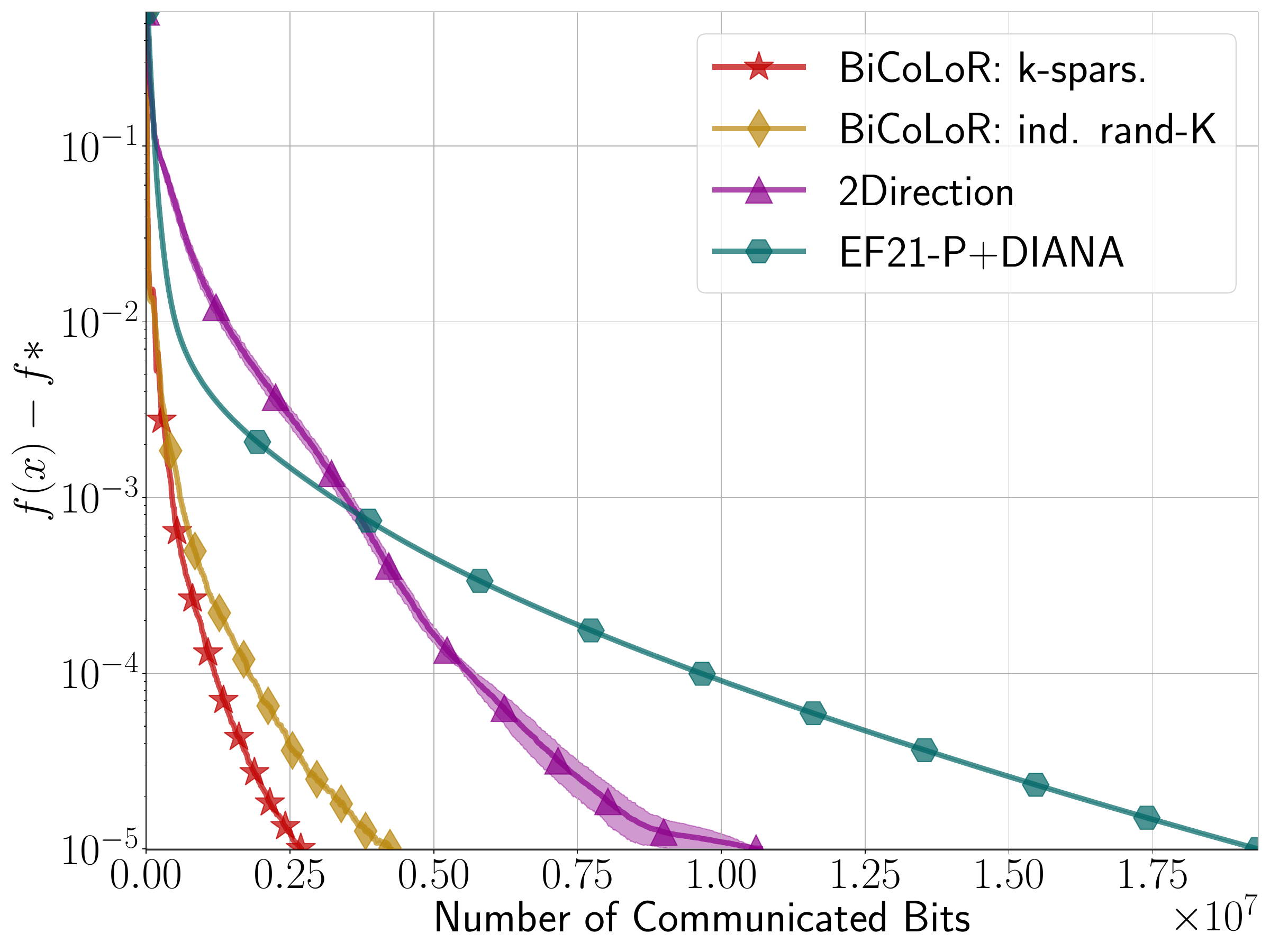}
    & 
    \includegraphics[width=0.47\linewidth]{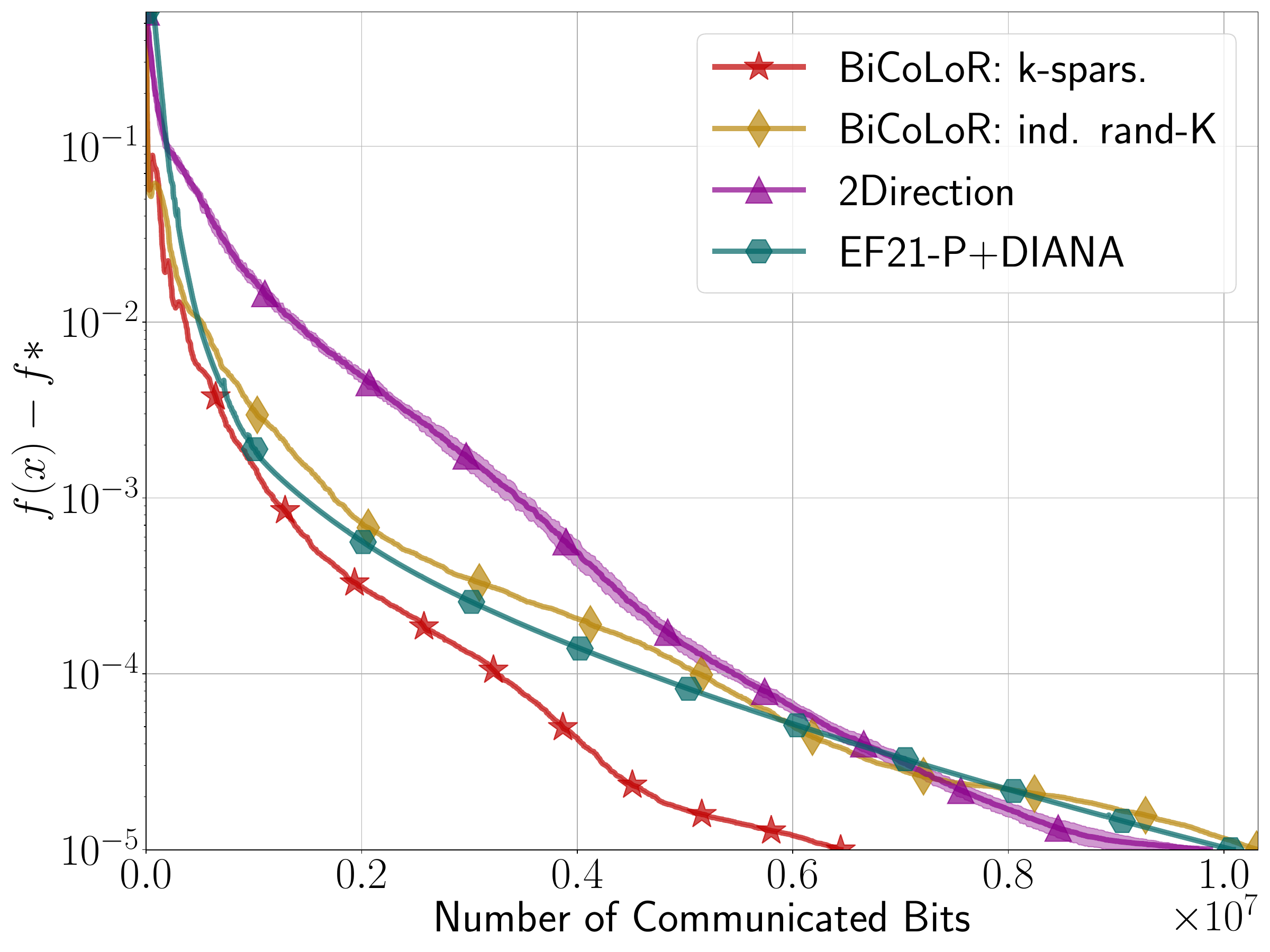} \\
    (a) $n=10$ &
    (b) $n=100$
    \end{tabular}
    \caption{
        Logistic regression on the \texttt{w8a} dataset.
        The compression scheme combines sparsification with $K = 100$ and Natural Compression.
    }
    \label{fig:w8a}
\end{figure}

Our experiments use datasets from the LibSVM library \citep{CC01a} 
(3-clause BSD license). Each dataset is first shuffled, 
then divided into $n$ pieces of same size $m$ assigned to the $n$ clients, discarding any leftover data points to ensure $m$ is an integer.

For \algno, we consider two compression strategies. We focus on the case $\alpha=1$, so that the compression level is the same downlink and uplink.
\begin{enumerate}
\item Sparsification is performed using a low value of the parameter $k$, like in Corollary~\ref{cor43}. The compressors are independent Natural compressors, performing quantization of reals on 9 bits, as explained in Section~\ref{secuc}.
\item $k=d$ and the compressors are the composition of independent \texttt{rand}-$K$ compression with same $K$, like in Corollary~\ref{cor42}, and Natural compression. See in Section~\ref{secuc} for the variance of a composition.
\end{enumerate}
For a given $K\in[d]$ and $k=K$ in the first strategy, the compression level is the same, but the difference is that every machine selects $K$ random coordinates independently in the second strategy, whereas the $K$ randomly selected coordinates are the same for all machines in the first one. Independence is beneficial but requires lower stepsizes $\eta$ and $\rho$, so we don't know a priori which of the two strategies is best.

For the other methods \algn{2Direction} and \algn{EF21-P+DIANA}, we also use a combination of \texttt{rand}-$K$ and Natural compression, both downlink and uplink 
(with appropriate scaling 
to make the downlink compressor contractive).

We tuned the stepsizes ($\gamma$ for \algno) for all algorithms. 
All other parameters are set to their best theoretical value. The algorithms are initialized with zero vectors.

Figures \ref{fig:real_sim} and \ref{fig:w8a} show the results on the  \texttt{real-sim} dataset (72,309 samples and $d=20$,958 features) and \texttt{w8a} dataset (49,749 samples and $d=300$ features), respectively.

\algno with the first strategy (using the parameter $k$ for sparsification instead of the compressors, red curves in the plots) outperforms the other algorithms. So, our theoretical findings are confirmed in practice, and \algno establishes a new state of the art for optimization with bidirectional compression.

\begin{figure}
    \centering
    \includegraphics[width=0.5\linewidth]{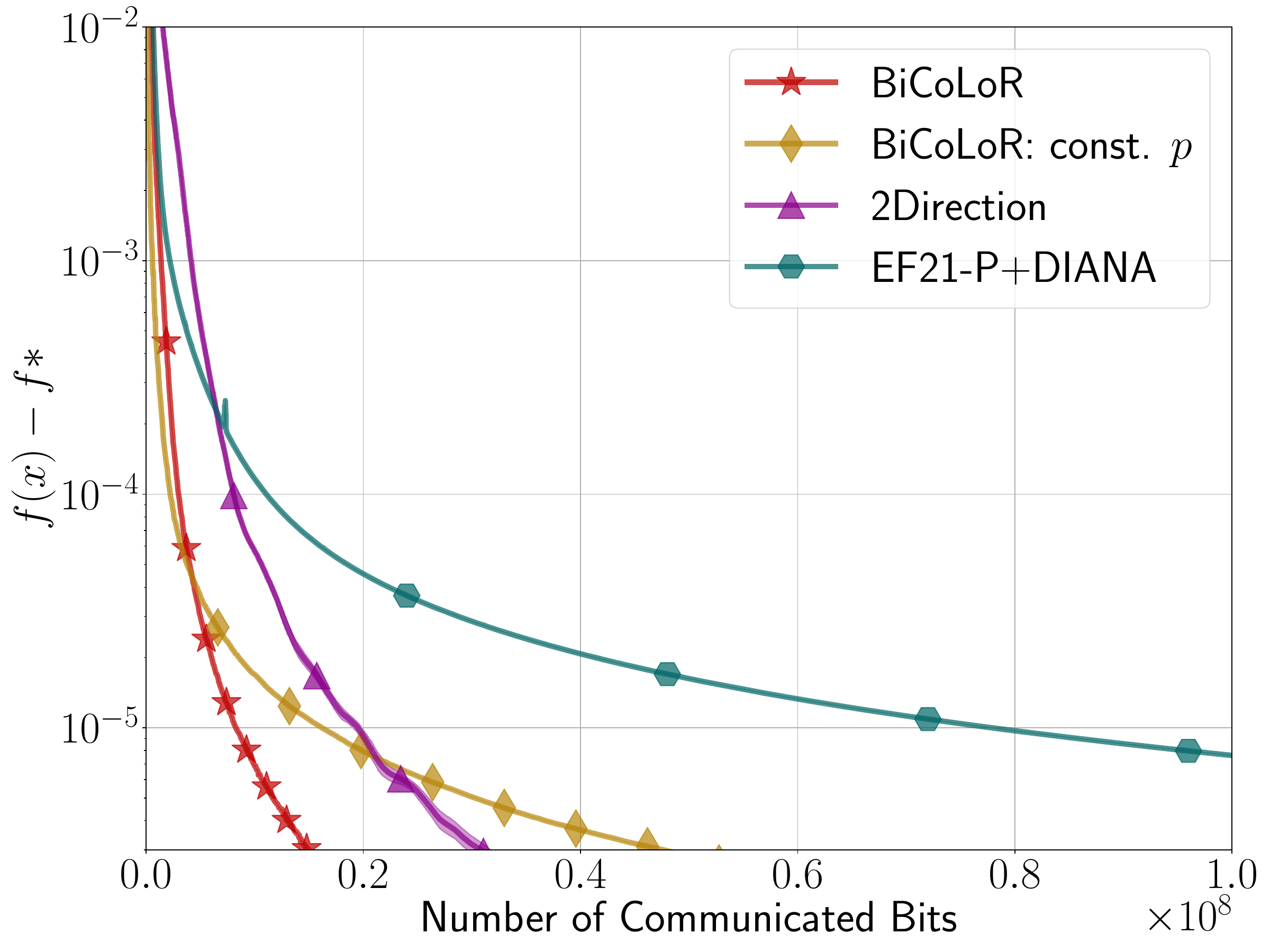}
    \caption{
        Logistic regression without regularization on the \texttt{w8a} dataset.
        The compression scheme combines sparsification with $K = 100$ and Natural Compression.
        \algno with decreasing $p_t$ (as defined in Theorem \ref{theogc}) outperforms the constant-$p$ variant in the long run.
    }
    \label{fig:convex}
\end{figure}

\subsection{General convex case}

In this section, we aim to demonstrate the acceleration benefit of using a decreasing sequence $p_t$, as defined in Theorem \ref{theogc}, in the general convex setting.
To this end, we conduct experiments using a convex loss function by removing the regularization term from the logistic regression objective \eqref{eqlogr}, i.e., setting $\mu = 0$.
The experiments are performed on the \texttt{w8a} dataset with $n = 10$.

We compare \algno with decreasing $p_t$ (as in Theorem \ref{theogc}) against \algno with constant $p_t \equiv p$, where $p$ is chosen so that both variants perform the same total number of local steps.
For compression, like in Section~\ref{secssc}, we consider sparsification with the parameter $k$ and Natural compression in \algno, and 
the other algorithms use 
compositions of \texttt{rand}-$K$ and Natural compression.
We tuned the stepsizes ($\gamma$ for \algno) for all methods, and all other parameters are set to their best theoretical values.

The results are shown in Figure \ref{fig:convex}. We observe that \algno with decreasing $p_t$ outperforms the constant-$p$ version, as well as the other algorithms.
Here too, our theoretical findings demonstrating acceleration are confirmed in practice, and \algno sets a new standard.

\section*{Acknowledgements}
This work was supported by funding from King Abdullah University of Science and Technology (KAUST):\\
i) KAUST Baseline Research Scheme,\\
ii) Center of Excellence for Generative AI (award no.\ 5940),\\
iii) Competitive Research Grant (CRG) Program (award no.\ 6460),\\
iv) SDAIA-KAUST Center of Excellence in Data Science and Artificial Intelligence (SDAIA-KAUST AI).

\bibliography{IEEEabrv,bib}

\begin{thebibliography}{63}
\providecommand{\natexlab}[1]{#1}
\providecommand{\url}[1]{\texttt{#1}}
\expandafter\ifx\csname urlstyle\endcsname\relax
  \providecommand{\doi}[1]{doi: #1}\else
  \providecommand{\doi}{doi: \begingroup \urlstyle{rm}\Url}\fi

\bibitem[Albasyoni et~al.(2020)Albasyoni, Safaryan, Condat, and
  {Richt\'arik}]{alb20}
Albasyoni, A., Safaryan, M., Condat, L., and {Richt\'arik}, P.
\newblock Optimal gradient compression for distributed and federated learning.
\newblock preprint arXiv:2010.03246, 2020.

\bibitem[Basu et~al.(2020)Basu, Data, Karakus, and Diggavi]{bas20}
Basu, D., Data, D., Karakus, C., and Diggavi, S.~N.
\newblock {Qsparse-Local-SGD: Distributed SGD With Quantization,
  Sparsification, and Local Computations}.
\newblock \emph{IEEE Journal on Selected Areas in Information Theory},
  1\penalty0 (1):\penalty0 217--226, 2020.

\bibitem[Bertsekas(2015)]{ber15}
Bertsekas, D.~P.
\newblock \emph{Convex optimization algorithms}.
\newblock Athena Scientific, Belmont, MA, USA, 2015.

\bibitem[Beznosikov et~al.(2020)Beznosikov, Horv{\'a}th, Richt{\'a}rik, and
  Safaryan]{bez20}
Beznosikov, A., Horv{\'a}th, S., Richt{\'a}rik, P., and Safaryan, M.
\newblock On biased compression for distributed learning.
\newblock preprint arXiv:2002.12410, 2020.

\bibitem[Bonawitz et~al.(2017)Bonawitz, Ivanov, Kreuter, Marcedone, {McMahan},
  Patel, Ramage, Segal, and Seth]{bon17}
Bonawitz, K., Ivanov, V., Kreuter, B., Marcedone, A., {McMahan}, H.~B., Patel,
  S., Ramage, D., Segal, A., and Seth, K.
\newblock Practical secure aggregation for privacy-preserving machine learning.
\newblock In \emph{Proc. of the 2017 ACM SIGSAC Conference on Computer and
  Communications Security}, pp.\  1175--1191, 2017.

\bibitem[Chang \& Lin(2011)Chang and Lin]{CC01a}
Chang, C.-C. and Lin, C.-J.
\newblock {LIBSVM}: {A} library for support vector machines.
\newblock \emph{ACM Transactions on Intelligent Systems and Technology},
  2:\penalty0 27:1--27:27, 2011.
\newblock Software available at http://www.csie.ntu.edu.tw/\%7Ecjlin/libsvm.

\bibitem[Condat \& {Richt\'arik}(2022)Condat and {Richt\'arik}]{con22mu}
Condat, L. and {Richt\'arik}, P.
\newblock {MURANA: A} generic framework for stochastic variance-reduced
  optimization.
\newblock In \emph{Proc. of the conference Mathematical and Scientific Machine
  Learning (MSML), PMLR 190}, 2022.

\bibitem[Condat \& Richt{\'a}rik(2023)Condat and Richt{\'a}rik]{con22rp}
Condat, L. and Richt{\'a}rik, P.
\newblock {RandProx}: {P}rimal-dual optimization algorithms with randomized
  proximal updates.
\newblock In \emph{Proc. of International Conference on Learning
  Representations (ICLR)}, 2023.

\bibitem[Condat et~al.(2022{\natexlab{a}})Condat, Agarsk{\'y}, and
  Richt{\'a}rik]{con22cs}
Condat, L., Agarsk{\'y}, I., and Richt{\'a}rik, P.
\newblock Provably doubly accelerated federated learning: {The} first
  theoretically successful combination of local training and compressed
  communication.
\newblock preprint arXiv:2210.13277, 2022{\natexlab{a}}.

\bibitem[Condat et~al.(2022{\natexlab{b}})Condat, Li, and
  {Richt\'arik}]{con22e}
Condat, L., Li, K., and {Richt\'arik}, P.
\newblock {EF-BV: A} unified theory of error feedback and variance reduction
  mechanisms for biased and unbiased compression in distributed optimization.
\newblock In \emph{Proc. of Conf. Neural Information Processing Systems
  (NeurIPS)}, 2022{\natexlab{b}}.

\bibitem[Condat et~al.(2023)Condat, Agarsk{\'y}, Malinovsky, and
  Richt{\'a}rik]{tamu23}
Condat, L., Agarsk{\'y}, I., Malinovsky, G., and Richt{\'a}rik, P.
\newblock {TAMUNA: Doubly} accelerated federated learning with local training,
  compression, and partial participation.
\newblock preprint arXiv:2302.09832 presented at the {\textit{Int. Workshop on
  Federated Learning in the Age of Foundation Models in Conjunction with
  NeurIPS 2023}}, 2023.

\bibitem[Condat et~al.(2025)Condat, Maranjyan, and Richt{\'a}rik]{con25loco}
Condat, L., Maranjyan, A., and Richt{\'a}rik, P.
\newblock {LoCoDL: C}ommunication-efficient distributed learning with local
  training and compression.
\newblock In \emph{Proc. of International Conference on Learning
  Representations (ICLR)}, 2025.

\bibitem[Egger et~al.(2025)Egger, Bitar, {Wachter-Zeh}, Weinberger, and
  {G\"und\"uz}]{egg25}
Egger, M., Bitar, R., {Wachter-Zeh}, A., Weinberger, N., and {G\"und\"uz}, D.
\newblock {BiCompFL: Stochastic} federated learning with bi-directional
  compression.
\newblock preprint arXiv:2502.00206, 2025.

\bibitem[Fatkhullin et~al.(2021)Fatkhullin, Sokolov, Gorbunov, Li, and
  Richt{\'a}rik]{fat21}
Fatkhullin, I., Sokolov, I., Gorbunov, E., Li, Z., and Richt{\'a}rik, P.
\newblock {EF21} with bells {\&} whistles: Practical algorithmic extensions of
  modern error feedback.
\newblock preprint arXiv:2110.03294, 2021.

\bibitem[Glasgow et~al.(2022)Glasgow, Yuan, and Ma]{gla22}
Glasgow, M.~R., Yuan, H., and Ma, T.
\newblock Sharp bounds for federated averaging {(Local SGD)} and continuous
  perspective.
\newblock In \emph{Proc. of Int. Conf. Artificial Intelligence and Statistics
  (AISTATS), PMLR 151}, pp.\  9050--9090, 2022.

\bibitem[Gorbunov et~al.(2020)Gorbunov, Hanzely, and {Richt\'{a}rik}]{gor202}
Gorbunov, E., Hanzely, F., and {Richt\'{a}rik}, P.
\newblock A unified theory of {SGD: Variance} reduction, sampling, quantization
  and coordinate descent.
\newblock In \emph{Proc. of 23rd Int. Conf. Artificial Intelligence and
  Statistics (AISTATS), PMLR 108}, 2020.

\bibitem[Gorbunov et~al.(2021)Gorbunov, Hanzely, and {Richt\'{a}rik}]{gor21}
Gorbunov, E., Hanzely, F., and {Richt\'{a}rik}, P.
\newblock {Local SGD: Unified} theory and new efficient methods.
\newblock In \emph{Proc. of 24th Int. Conf. Artificial Intelligence and
  Statistics (AISTATS), PMLR 130}, pp.\  3556--3564, 2021.

\bibitem[Gower et~al.(2020)Gower, Schmidt, Bach, and Richt\'{a}rik]{gow20a}
Gower, R.~M., Schmidt, M., Bach, F., and Richt\'{a}rik, P.
\newblock Variance-reduced methods for machine learning.
\newblock \emph{Proc. of the IEEE}, 108\penalty0 (11):\penalty0 1968--1983,
  November 2020.

\bibitem[Gruntkowska et~al.(2023)Gruntkowska, Tyurin, and Richt{á}rik]{gru23}
Gruntkowska, K., Tyurin, A., and Richt{á}rik, P.
\newblock {EF21-P} and friends: {I}mproved theoretical communication complexity
  for distributed optimization with bidirectional compression.
\newblock In \emph{Proc. of 40th Int. Conf. Machine Learning (ICML)}, 2023.

\bibitem[Gruntkowska et~al.(2024)Gruntkowska, Tyurin, and Richt{\'a}rik]{gru24}
Gruntkowska, K., Tyurin, A., and Richt{\'a}rik, P.
\newblock Improving the worst-case bidirectional communication complexity for
  nonconvex distributed optimization under function similarity.
\newblock In \emph{Proc. of Conf. Neural Information Processing Systems
  (NeurIPS)}, 2024.

\bibitem[Haddadpour \& Mahdavi(2019)Haddadpour and Mahdavi]{had19}
Haddadpour, F. and Mahdavi, M.
\newblock {On the Convergence of Local Descent Methods in Federated Learning}.
\newblock preprint arXiv:1910.14425, 2019.

\bibitem[Haddadpour et~al.(2021)Haddadpour, Kamani, Mokhtari, and
  Mahdavi]{had21}
Haddadpour, F., Kamani, M.~M., Mokhtari, A., and Mahdavi, M.
\newblock Federated learning with compression: Unified analysis and sharp
  guarantees.
\newblock In \emph{Proc. of Int. Conf. Artificial Intelligence and Statistics
  (AISTATS), PMLR 130}, pp.\  2350--2358, 2021.

\bibitem[Hanzely \& Richt{\'a}rik(2019)Hanzely and Richt{\'a}rik]{han19}
Hanzely, F. and Richt{\'a}rik, P.
\newblock One method to rule them all: Variance reduction for data, parameters
  and many new methods.
\newblock preprint arXiv:1905.11266, 2019.

\bibitem[He et~al.(2023{\natexlab{a}})He, Huang, Chen, Yin, and Yuan]{he23l}
He, Y., Huang, X., Chen, Y., Yin, W., and Yuan, K.
\newblock Lower bounds and accelerated algorithms in distributed stochastic
  optimization with communication compression.
\newblock preprint arXiv:2305.07612, 2023{\natexlab{a}}.

\bibitem[He et~al.(2023{\natexlab{b}})He, Huang, and Yuan]{yhe23}
He, Y., Huang, X., and Yuan, K.
\newblock Unbiased compression saves communication in distributed optimization:
  {W}hen and how much?
\newblock In \emph{Proc. of Conf. Neural Information Processing Systems
  (NeurIPS)}, 2023{\natexlab{b}}.

\bibitem[Horv\'{a}th et~al.(2022)Horv\'{a}th, Ho, Horv\'{a}th, Sahu, Canini,
  and Richt\'{a}rik]{hor22}
Horv\'{a}th, S., Ho, C.-Y., Horv\'{a}th, L., Sahu, A.~N., Canini, M., and
  Richt\'{a}rik, P.
\newblock Natural compression for distributed deep learning.
\newblock In \emph{Proc. of the conference Mathematical and Scientific Machine
  Learning (MSML), PMLR 190}, 2022.

\bibitem[{Horv\'ath} et~al.(2022){Horv\'ath}, Kovalev, Mishchenko, Stich, and
  {Richt\'arik}]{hor22d}
{Horv\'ath}, S., Kovalev, D., Mishchenko, K., Stich, S., and {Richt\'arik}, P.
\newblock Stochastic distributed learning with gradient quantization and
  variance reduction.
\newblock \emph{Optimization Methods and Software}, 2022.

\bibitem[Huang et~al.(2022)Huang, Chen, Yin, and Yuan]{hua22}
Huang, X., Chen, Y., Yin, W., and Yuan, K.
\newblock Lower bounds and nearly optimal algorithms in distributed learning
  with communication compression.
\newblock In \emph{Proc. of Conf. Neural Information Processing Systems
  (NeurIPS)}, pp.\  18955--18969, 2022.

\bibitem[Kairouz et~al.(2021)]{kai19}
Kairouz, P. et~al.
\newblock Advances and open problems in federated learning.
\newblock \emph{Foundations and Trends in Machine Learning}, 14\penalty0
  (1--2), 2021.

\bibitem[Karimireddy et~al.(2020)Karimireddy, Kale, Mohri, Reddi, Stich, and
  Suresh]{kar20}
Karimireddy, S.~P., Kale, S., Mohri, M., Reddi, S., Stich, S.~U., and Suresh,
  A.~T.
\newblock {SCAFFOLD}: {Stochastic} controlled averaging for federated learning.
\newblock In \emph{Proc. of 37th Int. Conf. Machine Learning (ICML)}, pp.\
  5132--5143, 2020.

\bibitem[Khaled et~al.(2019)Khaled, Mishchenko, and Richt\'{a}rik]{kha19}
Khaled, A., Mishchenko, K., and Richt\'{a}rik, P.
\newblock Better communication complexity for local {SGD}.
\newblock In \emph{NeurIPS Workshop on Federated Learning for Data Privacy and
  Confidentiality}, 2019.

\bibitem[Khaled et~al.(2020)Khaled, Mishchenko, and {Richt\'{a}rik}]{kha20a}
Khaled, A., Mishchenko, K., and {Richt\'{a}rik}, P.
\newblock Tighter theory for local {SGD} on identical and heterogeneous data.
\newblock In \emph{Proc. of 23rd Int. Conf. Artificial Intelligence and
  Statistics (AISTATS), PMLR 108}, 2020.

\bibitem[{Kone{\v{c}}n{\'y}} et~al.(2016{\natexlab{a}}){Kone{\v{c}}n{\'y}},
  {McMahan}, Ramage, and Richt{\'a}rik]{kon16a}
{Kone{\v{c}}n{\'y}}, J., {McMahan}, H.~B., Ramage, D., and Richt{\'a}rik, P.
\newblock Federated optimization: distributed machine learning for on-device
  intelligence.
\newblock arXiv:1610.02527, 2016{\natexlab{a}}.

\bibitem[{Kone{\v{c}}n{\'y}} et~al.(2016{\natexlab{b}}){Kone{\v{c}}n{\'y}},
  {McMahan}, Yu, Richt{\'a}rik, Suresh, and Bacon]{kon16}
{Kone{\v{c}}n{\'y}}, J., {McMahan}, H.~B., Yu, F.~X., Richt{\'a}rik, P.,
  Suresh, A.~T., and Bacon, D.
\newblock Federated learning: Strategies for improving communication
  efficiency.
\newblock In \emph{NIPS Private Multi-Party Machine Learning Workshop},
  2016{\natexlab{b}}.
\newblock arXiv:1610.05492.

\bibitem[Li et~al.(2020{\natexlab{a}})Li, Sahu, Talwalkar, and Smith]{li20}
Li, T., Sahu, A.~K., Talwalkar, A., and Smith, V.
\newblock Federated learning: Challenges, methods, and future directions.
\newblock \emph{IEEE Signal Processing Magazine}, 3\penalty0 (37):\penalty0
  50--60, 2020{\natexlab{a}}.

\bibitem[Li et~al.(2020{\natexlab{b}})Li, Huang, Yang, Wang, and Zhang]{li20a}
Li, X., Huang, K., Yang, W., Wang, S., and Zhang, Z.
\newblock On the convergence of {FedAvg} on non-iid data.
\newblock In \emph{Proc. of Int. Conf. Learning Representations (ICLR)},
  2020{\natexlab{b}}.

\bibitem[Li et~al.(2020{\natexlab{c}})Li, Kovalev, Qian, and
  Richt{á}rik]{zli20}
Li, Z., Kovalev, D., Qian, X., and Richt{á}rik, P.
\newblock Acceleration for compressed gradient descent in distributed and
  federated optimization.
\newblock In \emph{Proc. of 37th Int. Conf. Machine Learning (ICML)}, volume
  PMLR 119, 2020{\natexlab{c}}.

\bibitem[Liu et~al.(2020)Liu, Li, Tang, and Yan]{liu20}
Liu, X., Li, Y., Tang, J., and Yan, M.
\newblock A double residual compression algorithm for efficient distributed
  learning.
\newblock In \emph{Proc. of Int. Conf. Artificial Intelligence and Statistics
  (AISTATS), PMLR 108}, pp.\  133--143, 2020.

\bibitem[Malinovsky \& Richt{\'a}rik(2022)Malinovsky and
  Richt{\'a}rik]{malinovsky2022federated}
Malinovsky, G. and Richt{\'a}rik, P.
\newblock Federated random reshuffling with compression and variance reduction.
\newblock preprint arXiv:arXiv:2205.03914, 2022.

\bibitem[Malinovsky et~al.(2020)Malinovsky, Kovalev, Gasanov, Condat, and
  Richt{\'a}rik]{mal20}
Malinovsky, G., Kovalev, D., Gasanov, E., Condat, L., and Richt{\'a}rik, P.
\newblock From local {SGD} to local fixed point methods for federated learning.
\newblock In \emph{Proc. of 37th Int. Conf. Machine Learning (ICML)}, 2020.

\bibitem[Malinovsky et~al.(2022)Malinovsky, Yi, and {Richt\'arik}]{mal22}
Malinovsky, G., Yi, K., and {Richt\'arik}, P.
\newblock Variance reduced {ProxSkip}: {A}lgorithm, theory and application to
  federated learning.
\newblock In \emph{Proc. of Conf. Neural Information Processing Systems
  (NeurIPS)}, 2022.

\bibitem[Maranjyan et~al.(2022)Maranjyan, Safaryan, and
  Richtárik]{maranjyan2023gradskip}
Maranjyan, A., Safaryan, M., and Richtárik, P.
\newblock {GradSkip}: {C}ommunication-accelerated local gradient methods with
  better computational complexity.
\newblock preprint arXiv:2210.16402, 2022.

\bibitem[McMahan et~al.(2017)McMahan, Moore, Ramage, Hampson, and
  y~Arcas]{mcm17}
McMahan, H.~B., Moore, E., Ramage, D., Hampson, S., and y~Arcas, B.~A.
\newblock Communication-efficient learning of deep networks from decentralized
  data.
\newblock In \emph{Proc. of Int. Conf. Artificial Intelligence and Statistics
  (AISTATS), PMLR 54}, 2017.

\bibitem[Mishchenko et~al.(2022)Mishchenko, Malinovsky, Stich, and
  {Richt\'arik}]{mis22}
Mishchenko, K., Malinovsky, G., Stich, S., and {Richt\'arik}, P.
\newblock {ProxSkip: Yes! Local Gradient Steps Provably Lead to Communication
  Acceleration! Finally!}
\newblock In \emph{Proc. of the 39th International Conference on Machine
  Learning (ICML)}, July 2022.

\bibitem[Mishchenko et~al.(2024)Mishchenko, Gorbunov, {Tak\'a\v{c}}, and
  {Richt\'arik}]{mis19d}
Mishchenko, K., Gorbunov, E., {Tak\'a\v{c}}, M., and {Richt\'arik}, P.
\newblock Distributed learning with compressed gradient differences.
\newblock \emph{Optimization Methods and Software}, 2024.
\newblock Written in 2019 (arXiv:1901.09269).

\bibitem[Mitra et~al.(2021)Mitra, Jaafar, Pappas, and Hassani]{mit21}
Mitra, A., Jaafar, R., Pappas, G., and Hassani, H.
\newblock Linear convergence in federated learning: Tackling client
  heterogeneity and sparse gradients.
\newblock In \emph{Proc. of Conf. Neural Information Processing Systems
  (NeurIPS)}, 2021.

\bibitem[Philippenko \& Dieuleveut(2020)Philippenko and Dieuleveut]{phi20}
Philippenko, C. and Dieuleveut, A.
\newblock Bidirectional compression in heterogeneous settings for distributed
  or federated learning with partial participation: tight convergence
  guarantees.
\newblock preprint arXiv:2006.14591, 2020.

\bibitem[Philippenko \& Dieuleveut(2021)Philippenko and Dieuleveut]{phi21}
Philippenko, C. and Dieuleveut, A.
\newblock Preserved central model for faster bidirectional compression in
  distributed settings.
\newblock In \emph{Proc. of Conf. Neural Information Processing Systems
  (NeurIPS)}, 2021.

\bibitem[Reisizadeh et~al.(2020)Reisizadeh, Mokhtari, Hassani, Jadbabaie, and
  Pedarsani]{rei20}
Reisizadeh, A., Mokhtari, A., Hassani, H., Jadbabaie, A., and Pedarsani, R.
\newblock {FedPAQ}: {A} communication-efficient federated learning method with
  periodic averaging and quantization.
\newblock In \emph{Proc. of Int. Conf. Artificial Intelligence and Statistics
  (AISTATS)}, pp.\  2021--2031, 2020.

\bibitem[Richt{\'a}rik et~al.(2021)Richt{\'a}rik, Sokolov, and
  Fatkhullin]{ric21}
Richt{\'a}rik, P., Sokolov, I., and Fatkhullin, I.
\newblock {EF21}: {A} new, simpler, theoretically better, and practically
  faster error feedback.
\newblock In \emph{Proc. of 35th Conf. Neural Information Processing Systems
  (NeurIPS)}, 2021.

\bibitem[Sadiev et~al.(2022)Sadiev, Malinovsky, Gorbunov, Sokolov, Khaled,
  Burlachenko, and Richt\'arik]{sad22}
Sadiev, A., Malinovsky, G., Gorbunov, E., Sokolov, I., Khaled, A., Burlachenko,
  K., and Richt\'arik, P.
\newblock Federated optimization algorithms with random reshuffling and
  gradient compression.
\newblock preprint arXiv:2206.07021, 2022.

\bibitem[Safaryan et~al.(2022)Safaryan, Shulgin, and Richtárik]{saf22}
Safaryan, M., Shulgin, E., and Richtárik, P.
\newblock Uncertainty principle for communication compression in distributed
  and federated learning and the search for an optimal compressor.
\newblock \emph{Information and Inference: A Journal of the IMA}, 11\penalty0
  (2):\penalty0 557--580, 2022.

\bibitem[Scaman et~al.(2019)Scaman, Bach, Bubeck, Lee, and
  {Massouli\'e}]{sca19}
Scaman, K., Bach, F., Bubeck, S., Lee, Y.~T., and {Massouli\'e}, L.
\newblock Optimal convergence rates for convex distributed optimization in
  networks.
\newblock \emph{Journal of Machine Learning Research}, 20:\penalty0 1--31,
  2019.

\bibitem[Shalev-Shwartz \& Ben-David(2014)Shalev-Shwartz and
  Ben-David]{shai_book}
Shalev-Shwartz, S. and Ben-David, S.
\newblock \emph{Understanding machine learning: {F}rom theory to algorithms}.
\newblock Cambridge University Press, 2014.

\bibitem[Sra et~al.(2011)Sra, Nowozin, and Wright]{sra11}
Sra, S., Nowozin, S., and Wright, S.~J.
\newblock \emph{Optimization for Machine Learning}.
\newblock The MIT Press, 2011.

\bibitem[Stich(2019)]{sti19}
Stich, S.~U.
\newblock {Local SGD} converges fast and communicates little.
\newblock In \emph{Proc. of International Conference on Learning
  Representations (ICLR)}, 2019.

\bibitem[Tyurin \& Richt{á}rik(2023)Tyurin and Richt{á}rik]{tyu23}
Tyurin, A. and Richt{á}rik, P.
\newblock {2Direction}: {T}heoretically faster distributed training with
  bidirectional communication compression.
\newblock In \emph{Proc. of Conf. Neural Information Processing Systems
  (NeurIPS)}, 2023.

\bibitem[Wang et~al.(2021)]{wan21}
Wang, J. et~al.
\newblock A field guide to federated optimization.
\newblock preprint arXiv:2107.06917, 2021.

\bibitem[Woodworth et~al.(2020)Woodworth, Patel, and Srebro]{woo20}
Woodworth, B.~E., Patel, K.~K., and Srebro, N.
\newblock Minibatch vs {Local SGD} for heterogeneous distributed learning.
\newblock In \emph{Proc. of Conf. Neural Information Processing Systems
  (NeurIPS)}, 2020.

\bibitem[Yang et~al.(2021)Yang, Fang, and Liu]{yang21}
Yang, H., Fang, M., and Liu, J.
\newblock Achieving linear speedup with partial worker participation in
  non-{IID} federated learning.
\newblock In \emph{Proc. of International Conference on Learning
  Representations (ICLR)}, 2021.

\bibitem[Yang et~al.(2024)Yang, Qiu, Khanduri, Fang, and Liu]{yang23}
Yang, H., Qiu, P., Khanduri, P., Fang, M., and Liu, J.
\newblock Understanding server-assisted federated learning in the presence of
  incomplete client participation.
\newblock In \emph{Proc. of Int. Conf. Machine Learning (ICML)}, 2024.

\bibitem[Yi et~al.(2025)Yi, Condat, and Richt{\'a}rik]{scafflix}
Yi, K., Condat, L., and Richt{\'a}rik, P.
\newblock Explicit personalization and local training: {D}ouble communication
  acceleration in federated learning.
\newblock \emph{Transactions on Machine Learning Research}, June 2025.

\bibitem[Zhao et~al.(2018)Zhao, Li, Lai, Suda, Civin, and Chandra]{zhao18}
Zhao, Y., Li, M., Lai, L., Suda, N., Civin, D., and Chandra, V.
\newblock Federated learning with non-iid data.
\newblock preprint arXiv:1806.00582, 2018.

\end{thebibliography}
\bibliographystyle{icml2024}

\clearpage
\appendix

{\noindent\huge \textbf{Appendix}}

\small

\section{Proof of Theorem~\ref{theo1}}\label{secalg2}

We define the Euclidean space $\boldsymbol{\mathcal{X}}\eqdef {(\mathbb{R}^d)}^{n+2}$ endowed with the weighted inner product
\begin{equation*}
\langle \mathbf{x},\mathbf{x}'\rangle_{\boldsymbol{\mathcal{X}}}\eqdef  \sum_{i=1}^n \langle x_i,x_i'\rangle + 2n \langle x_s,x_s'\rangle
+ n \langle y,y'\rangle,
\end{equation*}
for every $\mathbf{x}=(x_1,\ldots,x_n,x_s,y)$, $\mathbf{x}'=(x_1',\ldots,x_n',x_s',y')$. We also define the Euclidean spaces $\boldsymbol{\mathcal{U}}\eqdef {(\mathbb{R}^d)}^{n}$ endowed with the standard inner product
$
\langle \mathbf{u},\mathbf{u}'\rangle_{\boldsymbol{\mathcal{U}}}\eqdef  \sum_{i=1}^n \langle u_i,u_i'\rangle 
$, 
for every $\mathbf{u}=(u_1,\ldots,u_n)$, $\mathbf{u}'=(u_1',\ldots,u_n')$, and $\boldsymbol{\mathcal{U}}_y\eqdef \mathbb{R}^d$ endowed with the weighted inner product
$
\langle u,u'\rangle_{\boldsymbol{\mathcal{U}}_y}\eqdef  n \langle u,u'\rangle 
$, for every $u,u'\in \mathbb{R}^d$.

We reformulate the problem  \eqref{eq1} as
\begin{equation}
\min_{\mathbf{x}\in \boldsymbol{\mathcal{X}}} \ \mathbf{f}(\mathbf{x})\quad\mbox{s.t.}\quad D\mathbf{x}=\mathbf{0}\ \mbox{ and }\ D_y\mathbf{x}=\mathbf{0},\label{eq2}
\end{equation}
where 
\begin{align*}
\mathbf{f}&:\mathbf{x}=(x_1,\ldots,x_n,x_s,y)\in \boldsymbol{\mathcal{X}} \mapsto \sum_{i=1}^n f_i(x_i) + 2nf_s(x_s)+ng(y),\\
D&:\boldsymbol{\mathcal{X}}\rightarrow \boldsymbol{\mathcal{U}} : (x_1,\ldots,x_n,x_s,y)\mapsto(x_1-x_s,\ldots,x_n-x_s),\\
D_y&:\boldsymbol{\mathcal{X}}\rightarrow \boldsymbol{\mathcal{U}}_y : (x_1,\ldots,x_n,x_s,y)\mapsto(y-x_s),
\end{align*}
and  $\mathbf{0}$ denotes the zero element of the Euclidean space.

We note that the function $\mathbf{f}$ in $\boldsymbol{\mathcal{X}}$  is $L$-smooth and $\mu$-strongly convex, and

\noindent$\nabla \mathbf{f}(\mathbf{x})=\big(\nabla f_1(x_1),\ldots \nabla f_n(x_n),\nabla f_s(x_s),\nabla g(y)\big)$. The adjoint operators are 
\begin{align*}
D^*&: \boldsymbol{\mathcal{U}} \rightarrow \boldsymbol{\mathcal{X}}: (u_1,\ldots,u_n)\mapsto\left(u_1,\ldots,u_n,-\frac{1}{2n}\sum_{i=1}^n u_i,0\right),\\
D_y^*&: \boldsymbol{\mathcal{U}}_y \rightarrow \boldsymbol{\mathcal{X}}: u_y\mapsto(0,\ldots,0,-\frac{1}{2}u_y,u_y).
\end{align*}

We also introduce vector notations for the variables in \algno. We define  $\mathbf{x}^\star\eqdef 
(x^\star,\ldots,x^\star,x^\star,x^\star)\in \boldsymbol{\mathcal{X}}$ as the unique solution to \eqref{eq2},
$\mathbf{w}^\star\eqdef \mathbf{x}^\star - \gamma \nabla \mathbf{f}(\mathbf{x}^\star)$
 and, 
for every $t\geq 0$, 
$\mathbf{x}^t\eqdef (x_1^t,\ldots,x_n^t,x_s^t,y^t)$, $\mathbf{\hat{x}}^t\eqdef (\hat{x}_1^t,\ldots,\hat{x}_n^t,\hat{x}_s^t,\hat{y}^t)$, $\mathbf{u}^t\eqdef (u_1^t,\ldots,u_n^t)$, $\mathbf{u}^\star\eqdef (u_1^\star,\ldots,u_n^\star)=\nabla \mathbf{f}(\mathbf{x}^\star)$, $\mathbf{w}^t\eqdef \mathbf{x}^t - \gamma \nabla \mathbf{f}(\mathbf{x}^t)$,
and the $\sigma$-algebra $\mathcal{F}^t$  generated by the collection of random variables $\mathbf{x}^0,\mathbf{u}^0,u_y^0\ldots, \mathbf{x}^t,\mathbf{u}^t,u_y^t$.

We first consider the case $k=d$, so that $\Omega^t=[d]$ and full vectors are compressed.

Let $t\geq 0$. 
We can write the iteration of \algno as 
\begin{equation}
\left\lfloor\begin{array}{l}
\mathbf{\hat{x}}^t\eqdef \mathbf{x}^t-\gamma \nabla \mathbf{f}(\mathbf{x}^t)+\gamma D^*\mathbf{u}^t+\gamma D_y^* u_y^t=\mathbf{w}^t+\gamma D^*\mathbf{u}^t+\gamma D_y^* u_y^t\\
\mbox{flip a coin $\theta^t\in\{0,1\}$, with $\mathrm{Prob}(\theta^t=1)=p$}\\
\mbox{\textbf{if} }\theta^t=1\\
\ \ \ \mathbf{r}^t\eqdef \big(\hat{x}_1^t-c_s^t-\hat{y}^t,\ldots,\hat{x}_n^t-c_s^t-\hat{y}^t,\frac{1}{2}(\hat{x}_s^t-\bar{c}^t-\hat{y}^t),0\big)\\
\ \ \ \mathbf{r}_y^t\eqdef \big(0,\ldots,0,\frac{1}{2}(\hat{x}_s^t-\hat{y}^t),-c_s^t\big)\\
\ \ \ \mathbf{x}^{t+1}\eqdef \mathbf{\hat{x}}^t-\rho\mathbf{r}
-\rho_y\mathbf{r}_y\\
\ \ \ \mathbf{u}^{t+1}\eqdef \mathbf{u}^{t}- \frac{p\eta}{\gamma}(c_1^t-c_s^t,\ldots,c_n^t-c_s^t)\\
\ \ \ u_y^{t+1}\eqdef u_y^{t}+ \frac{p\eta_y}{\gamma}c_s^t\\
\mbox{\textbf{else}}\\
\ \ \ \mathbf{x}^{t+1}\eqdef \mathbf{\hat{x}}^t\\
\ \ \ \mathbf{u}^{t+1}\eqdef \mathbf{u}^t, u_y^{t+1}=u_y^t\\
\mbox{\textbf{end if}}\\
\end{array}\right.
\end{equation}
We have
\begin{equation}
\Exp{\mathbf{x}^{t+1}\;|\;\mathcal{F}^t,\theta=1}=\mathbf{\hat{x}}^t-\rho D^*D\mathbf{\hat{x}}^t-\rho_y D_y^*D_y\mathbf{\hat{x}}^t.\label{equ1}
\end{equation}
That is, $\mathbf{x}^{t+1}$ is updated using a stochastic unbiased estimate of the right hand side of \eqref{equ1} made from the compressed vectors, because the server does not know the $\hat{x}_i^t$ and the clients do not know $\hat{x}_s^t$, the only available information consists of the compressed vectors. Similarly,
\begin{align*}
\Exp{\mathbf{u}^{t+1}\;|\;\mathcal{F}^t,\theta=1}&=\mathbf{u}^t- \frac{p\eta}{\gamma}D\mathbf{\hat{x}}^t\\
\Exp{u_y^{t+1}\;|\;\mathcal{F}^t,\theta=1}&=u_y^t - \frac{p\eta_y}{\gamma}D_y\mathbf{\hat{x}}^t.
\end{align*}
Thus,
\begin{align*}
\Expc{\sqn{\mathbf{x}^{t+1}-\mathbf{x}^\star}_{\boldsymbol{\mathcal{X}}}}&=(1-p)\sqn{\mathbf{\hat{x}}^t-\mathbf{x}^\star}_{\boldsymbol{\mathcal{X}}}+p\Exp{\sqn{\mathbf{x}^{t+1}-\mathbf{x}^\star}_{\boldsymbol{\mathcal{X}}}\;|\;\mathcal{F}^t,\theta^t=1}\\
&= (1-p)\sqn{\mathbf{\hat{x}}^t-\mathbf{x}^\star}_{\boldsymbol{\mathcal{X}}}+p\sqn{\mathbf{\hat{x}}^t-\mathbf{\hat{x}}^\star-\rho D^*D\mathbf{\hat{x}}^t-\rho_y D_y^*D_y\mathbf{\hat{x}}^t}_{\boldsymbol{\mathcal{X}}}\\
&\quad + pn(\rho^2+\rho_y^2) \Exp{\sqn{c_s^t-(\hat{x}_s^t-\hat{y}^t)}\;|\;\mathcal{F}^t,\theta^t=1}\\
&\quad + p\frac{2n\rho^2}{4} \Exp{\sqn{\bar{c}^t-\frac{1}{n}\sum_{i=1}^n (\hat{x}_i^t-\hat{y}^t)}\;|\;\mathcal{F}^t,\theta^t=1}\\
&\leq (1-p)\sqn{\mathbf{\hat{x}}^t-\mathbf{x}^\star}_{\boldsymbol{\mathcal{X}}}+p
\sqn{\mathbf{\hat{x}}^t-\mathbf{\hat{x}}^\star}_{\boldsymbol{\mathcal{X}}}\\
&\quad+p\sqn{\rho D^*D\mathbf{\hat{x}}^t+\rho_y D_y^*D_y\mathbf{\hat{x}}^t}_{\boldsymbol{\mathcal{X}}}\\
&\quad-2p\left\langle \mathbf{\hat{x}}^t-\mathbf{\hat{x}}^\star,\rho D^*D\mathbf{\hat{x}}^t+\rho_y D_y^*D_y\mathbf{\hat{x}}^t\right\rangle_{\boldsymbol{\mathcal{X}}}\\
&\quad + pn(\rho^2+\rho_y^2)\omega_s  \sqn{\hat{x}_s^t-\hat{y}^t}+ pn\rho^2  \frac{\oma }{2n} \sum_{i=1}^n \sqn{\hat{x}_i^t-\hat{y}^t}\\
&= \sqn{\mathbf{\hat{x}}^t-\mathbf{x}^\star}_{\boldsymbol{\mathcal{X}}}+p\sqn{\rho D^*D\mathbf{\hat{x}}^t+\rho_y D_y^*D_y\mathbf{\hat{x}}^t}_{\boldsymbol{\mathcal{X}}}\\
&\quad-2p\rho\sqn{D\mathbf{\hat{x}}^t}_{\boldsymbol{\mathcal{U}}}-2p\rho_y\sqn{D_y\mathbf{\hat{x}}^t}_{\boldsymbol{\mathcal{U}}_y}\\
&\quad + pn(\rho^2+\rho_y^2)\omega_s  \sqn{\hat{x}_s^t-\hat{y}^t} +  \frac{ p\rho^2\oma}{2} \sum_{i=1}^n \sqn{\hat{x}_i^t-\hat{y}^t}.
\end{align*}
Moreover,
\begin{align*}
 \sqn{\mathbf{\hat{x}}^t-\mathbf{x}^\star}_{\boldsymbol{\mathcal{X}}}
 &= \sqn{\mathbf{w}^t-\mathbf{w}^\star}_{\boldsymbol{\mathcal{X}}}+\gamma^2
 \sqn{D^*(\mathbf{u}^t-\mathbf{u}^\star)+D_y^*(u_y^t-u_y^\star)}_{\boldsymbol{\mathcal{X}}}\\
 &\quad+2\gamma\big\langle \mathbf{w}^t-\mathbf{w}^\star, D^*(\mathbf{u}^t-\mathbf{u}^\star)+D_y^*(u_y^t-u_y^\star)\big\rangle_{\boldsymbol{\mathcal{X}}}.
\end{align*}
On the other hand, 
\begin{align*}
\Expc{\sqn{\mathbf{u}^{t+1}-\mathbf{u}^\star}_{\boldsymbol{\mathcal{U}}}}
&=(1-p)\sqn{\mathbf{u}^t-\mathbf{u}^\star}_{\boldsymbol{\mathcal{U}}}+p\Exp{\sqn{\mathbf{u}^{t+1}-\mathbf{u}^\star}_{\boldsymbol{\mathcal{U}}}\;|\;\mathcal{F}^t,\theta^t=1}\\
&=(1-p)\sqn{\mathbf{u}^t-\mathbf{u}^\star}_{\boldsymbol{\mathcal{U}}}+p\sqn{\mathbf{u}^t-\mathbf{u}^\star- \frac{p\eta}{\gamma}D\mathbf{\hat{x}}^t}_{\boldsymbol{\mathcal{U}}}\\
&\quad+\frac{p^3\eta^2}{\gamma^2}\Exp{\sum_{i=1}^n\sqn{c_i^t-c_s^t-(\hat{x}_i^t-\hat{x}_s^t)}\;|\;\mathcal{F}^t,\theta^t=1}.
\end{align*}
Let $i\in [n]$. Since $\mathcal{C}_i$ and $\mathcal{C}_s$ are supposed independent, we have
\begin{align*}
\Exp{\sqn{c_i^t-c_s^t-(\hat{x}_i^t-\hat{x}_s^t)}\;|\;\mathcal{F}^t,\theta^t=1}&=\Exp{\sqn{c_i^t-(\hat{x}_i^t-\hat{y}^t)}\;|\;\mathcal{F}^t,\theta^t=1}\\
&\quad+\Exp{\sqn{c_s^t-(\hat{x}_s^t-\hat{y}^t)}\;|\;\mathcal{F}^t,\theta^t=1}\\
&\leq \omega \sqn{\hat{x}_i^t-\hat{y}^t}+\omega_s \sqn{\hat{x}_s^t-\hat{y}^t}.
\end{align*}
Therefore,
\begin{align*}
\Expc{\sqn{\mathbf{u}^{t+1}-\mathbf{u}^\star}_{\boldsymbol{\mathcal{U}}}}
&\leq(1-p)\sqn{\mathbf{u}^t-\mathbf{u}^\star}_{\boldsymbol{\mathcal{U}}}+p\sqn{\mathbf{u}^t-\mathbf{u}^\star}_{\boldsymbol{\mathcal{U}}}
+\frac{p^3\eta^2}{\gamma^2}\sqn{D\mathbf{\hat{x}}^t}_{\boldsymbol{\mathcal{U}}}\\
&\quad-\frac{2p^2\eta}{\gamma}\langle \mathbf{u}^t-\mathbf{u}^\star,D\mathbf{\hat{x}}^t\rangle_{\boldsymbol{\mathcal{U}}}\\
&\quad+\frac{p^3\eta^2\omega}{\gamma^2}\sum_{i=1}^n\sqn{\hat{x}_i^t-\hat{y}^t}+\frac{p^3\eta^2\omega_s n}{\gamma^2}\sqn{\hat{x}_s^t-\hat{y}^t}.
\end{align*}
Moreover,
\begin{align*}
\Expc{\sqn{u_y^{t+1}-u_y^\star}_{\boldsymbol{\mathcal{U}}_y}}
&=(1-p)\sqn{\mathbf{u}^t-\mathbf{u}^\star}_{\boldsymbol{\mathcal{U}}_y}+p\Exp{\sqn{\mathbf{u}^{t+1}-\mathbf{u}^\star}_{\boldsymbol{\mathcal{U}}_y}\;|\;\mathcal{F}^t,\theta^t=1}\\
&=(1-p)\sqn{u_y^t-u_y^\star}_{\boldsymbol{\mathcal{U}}_y}+p\sqn{u_y^t-u_y^\star- \frac{p\eta_y}{\gamma}D_y\mathbf{\hat{x}}^t}_{\boldsymbol{\mathcal{U}}_y}\\
&\quad+\frac{p^3\eta_y^2n}{\gamma^2}\Exp{\sqn{c_s^t-(\hat{x}_s^t-\hat{y}^t)}\;|\;\mathcal{F}^t,\theta^t=1}\\
&\leq (1-p)\sqn{u_y^t-u_y^\star}_{\boldsymbol{\mathcal{U}}_y}+p\sqn{u_y^t-u_y^\star}_{\boldsymbol{\mathcal{U}}_y}+\frac{p^3\eta_y^2}{\gamma^2}\sqn{D_y\mathbf{\hat{x}}^t}_{\boldsymbol{\mathcal{U}}_y}\\
&\quad-\frac{2p^2\eta_y}{\gamma}\langle u_y^t-u_y^\star,D_y\mathbf{\hat{x}}^t\rangle_{\boldsymbol{\mathcal{U}}_y}+\frac{p^3\eta_y^2\omega_s n}{\gamma^2}\sqn{\hat{x}_s^t-\hat{y}^t}.
\end{align*}
Thus, using the fact that $D\mathbf{x}^\star=D_y\mathbf{x}^\star=\mathbf{0}$, we have
\begin{align*}
&\frac{\gamma}{p^2\eta}\Expc{\sqn{\mathbf{u}^{t+1}-\mathbf{u}^\star}_{\boldsymbol{\mathcal{U}}}}+
\frac{\gamma}{p^2\eta_y}\Expc{\sqn{u_y^{t+1}-u_y^\star}_{\boldsymbol{\mathcal{U}}_y}}\\
&\leq\frac{\gamma}{p^2\eta}\sqn{\mathbf{u}^t-\mathbf{u}^\star}_{\boldsymbol{\mathcal{U}}}+\frac{\gamma}{p^2\eta_y}\sqn{u_y^t-u_y^\star}_{\boldsymbol{\mathcal{U}}_y}
-2\big\langle D^* (\mathbf{u}^t-\mathbf{u}^\star)+D_y^* (u_y^t-u_y^\star),\mathbf{\hat{x}}^t-\mathbf{x}^\star\big\rangle_{\boldsymbol{\mathcal{X}}}\\
&\quad+\frac{p\eta}{\gamma}\sqn{D\mathbf{\hat{x}}^t}_{\boldsymbol{\mathcal{U}}}+\frac{p\eta_y}{\gamma}\sqn{D_y\mathbf{\hat{x}}^t}_{\boldsymbol{\mathcal{U}}_y}+\frac{p\eta\omega}{\gamma}\sum_{i=1}^n\sqn{\hat{x}_i^t-\hat{y}^t}+\frac{p(\eta+\eta_y)\omega_s}{\gamma}\sqn{D_y\mathbf{\hat{x}}^t}_{\boldsymbol{\mathcal{U}}_y}\\
&=\frac{\gamma}{p^2\eta}\sqn{\mathbf{u}^t-\mathbf{u}^\star}_{\boldsymbol{\mathcal{U}}}+\frac{\gamma}{p^2\eta_y}\sqn{u_y^t-u_y^\star}_{\boldsymbol{\mathcal{U}}_y}-2\gamma \sqn{D^* (\mathbf{u}^t-\mathbf{u}^\star)+D_y^* (u_y^t-u_y^\star)}_{\boldsymbol{\mathcal{X}}}\\
&\quad-2\big\langle D^* (\mathbf{u}^t-\mathbf{u}^\star)+D_y^* (u_y^t-u_y^\star),\mathbf{w}^t-\mathbf{w}^\star\big\rangle_{\boldsymbol{\mathcal{X}}}\\
&\quad+\frac{p\eta}{\gamma}\sqn{D\mathbf{\hat{x}}^t}_{\boldsymbol{\mathcal{U}}}+\frac{p\eta_y}{\gamma}\sqn{D_y\mathbf{\hat{x}}^t}_{\boldsymbol{\mathcal{U}}_y}+\frac{p\eta\omega}{\gamma}\sum_{i=1}^n\sqn{\hat{x}_i^t-\hat{y}^t}+\frac{p(\eta+\eta_y)\omega_s}{\gamma}\sqn{D_y\mathbf{\hat{x}}^t}_{\boldsymbol{\mathcal{U}}_y}.
\end{align*}
Hence,
\begin{align*}
&\frac{1}{\gamma}\Expc{\sqn{\mathbf{x}^{t+1}-\mathbf{x}^\star}_{\boldsymbol{\mathcal{X}}}}+\frac{\gamma}{p^2\eta}\Expc{\sqn{\mathbf{u}^{t+1}-\mathbf{u}^\star}_{\boldsymbol{\mathcal{U}}}}+
\frac{\gamma}{p^2\eta_y}\Expc{\sqn{u_y^{t+1}-u_y^\star}_{\boldsymbol{\mathcal{U}}_y}}\\
&\leq \frac{1}{\gamma}\sqn{\mathbf{w}^t-\mathbf{w}^\star}_{\boldsymbol{\mathcal{X}}}+\gamma
 \sqn{D^*(\mathbf{u}^t-\mathbf{u}^\star)+D_y^*(u_y^t-u_y^\star)}_{\boldsymbol{\mathcal{X}}}\\
 &\quad+2\big\langle \mathbf{w}^t-\mathbf{w}^\star, D^*(\mathbf{u}^t-\mathbf{u}^\star)+D_y^*(u_y^t-u_y^\star)\big\rangle_{\boldsymbol{\mathcal{X}}}\\
&+\frac{\gamma}{p^2\eta}\sqn{\mathbf{u}^t-\mathbf{u}^\star}_{\boldsymbol{\mathcal{U}}}+\frac{\gamma}{p^2\eta_y}\sqn{u_y^t-u_y^\star}_{\boldsymbol{\mathcal{U}}_y}-2\gamma \sqn{D^* (\mathbf{u}^t-\mathbf{u}^\star)+D_y^* (u_y^t-u_y^\star)}_{\boldsymbol{\mathcal{X}}}\\
&\quad-2\big\langle D^* (\mathbf{u}^t-\mathbf{u}^\star)+D_y^* (u_y^t-u_y^\star),\mathbf{w}^t-\mathbf{w}^\star\big\rangle_{\boldsymbol{\mathcal{X}}}\\
&\quad+\frac{p}{\gamma}\sqn{\rho D^*D\mathbf{\hat{x}}^t+\rho_y D_y^*D_y\mathbf{\hat{x}}^t}_{\boldsymbol{\mathcal{X}}}
-\frac{2p\rho}{\gamma}\sqn{D\mathbf{\hat{x}}^t}_{\boldsymbol{\mathcal{U}}}-\frac{2p\rho_y}{\gamma}\sqn{D_y\mathbf{\hat{x}}^t}_{\boldsymbol{\mathcal{U}}_y}\\
&\quad + \frac{pn(\rho^2+\rho_y^2)\omega_s}{\gamma}  \sqn{\hat{x}_s^t-\hat{y}^t} + \frac{p\rho^2  \oma }{2\gamma} \sum_{i=1}^n \sqn{\hat{x}_i^t-\hat{y}^t}\\
&\quad+\frac{p\eta}{\gamma}\sqn{D\mathbf{\hat{x}}^t}_{\boldsymbol{\mathcal{U}}}+\frac{p\eta_y}{\gamma}\sqn{D_y\mathbf{\hat{x}}^t}_{\boldsymbol{\mathcal{U}}_y}+\frac{p\eta\omega}{\gamma}\sum_{i=1}^n\sqn{\hat{x}_i^t-\hat{y}^t}+\frac{p(\eta+\eta_y)\omega_s}{\gamma}\sqn{D_y\mathbf{\hat{x}}^t}_{\boldsymbol{\mathcal{U}}_y}\\
&= \frac{1}{\gamma}\sqn{\mathbf{w}^t-\mathbf{w}^\star}_{\boldsymbol{\mathcal{X}}}+\frac{\gamma}{p^2\eta}\sqn{\mathbf{u}^t-\mathbf{u}^\star}_{\boldsymbol{\mathcal{U}}}+\frac{\gamma}{p^2\eta_y}\sqn{u_y^t-u_y^\star}_{\boldsymbol{\mathcal{U}}_y}\\
&\quad-\gamma \sqn{D^* (\mathbf{u}^t-\mathbf{u}^\star)+D_y^* (u_y^t-u_y^\star)}_{\boldsymbol{\mathcal{X}}}+\frac{p}{\gamma}\sqn{\rho D^*D\mathbf{\hat{x}}^t+\rho_y D_y^*D_y\mathbf{\hat{x}}^t}_{\boldsymbol{\mathcal{X}}}\\
&\quad+\frac{p\eta-2p\rho}{\gamma}\sqn{D\mathbf{\hat{x}}^t}_{\boldsymbol{\mathcal{U}}}+ \frac{p\rho^2 \oma+2p\eta\omega}{2\gamma} \sum_{i=1}^n \sqn{\hat{x}_i^t-\hat{y}^t}\\
&\quad+\frac{p\eta_y-2p\rho_y+p(\eta+\eta_y)\omega_s+p(\rho^2+\rho_y^2)\omega_s}{\gamma}\sqn{D_y\mathbf{\hat{x}}^t}_{\boldsymbol{\mathcal{U}}_y}.\end{align*}
Using Young's inequality, we have
\begin{align*}
\sum_{i=1}^n \sqn{\hat{x}_i^t-\hat{y}^t}&\leq \sum_{i=1}^n \left(2\sqn{\hat{x}_i^t-\hat{x}_s^t}+2\sqn{\hat{x}_s^t-\hat{y}^t}\right)\\
&=2\sqn{D\mathbf{\hat{x}}^t}_{\boldsymbol{\mathcal{U}}}+2\sqn{D_y\mathbf{\hat{x}}^t}_{\boldsymbol{\mathcal{U}}_y}.
\end{align*}
The leading eigenvalue of $D^*D+ D_y^*D_y$ is 2, with corresponding eigenvector $(x,\ldots,x,-x,x)$ for any $x\in\mathbb{R}^d$, so that
\begin{align*}
\sqn{\rho D^*D\mathbf{\hat{x}}^t+\rho_y D_y^*D_y\mathbf{\hat{x}}^t}_{\boldsymbol{\mathcal{X}}}&\leq 
\|D^*D+ D_y^*D_y\|\left(\sqn{\rho D\mathbf{\hat{x}}^t}_{\boldsymbol{\mathcal{U}}}+\sqn{\rho_y D_y\mathbf{\hat{x}}^t}_{\boldsymbol{\mathcal{U}}_y}\right)\\
&=2\rho^2\sqn{D\mathbf{\hat{x}}^t}_{\boldsymbol{\mathcal{U}}}+2\rho_y^2\sqn{D_y\mathbf{\hat{x}}^t}_{\boldsymbol{\mathcal{U}}_y}.
\end{align*}
Therefore,
\begin{align*}
&\frac{1}{\gamma}\Expc{\sqn{\mathbf{x}^{t+1}-\mathbf{x}^\star}_{\boldsymbol{\mathcal{X}}}}+\frac{\gamma}{p^2\eta}\Expc{\sqn{\mathbf{u}^{t+1}-\mathbf{u}^\star}_{\boldsymbol{\mathcal{U}}}}+
\frac{\gamma}{p^2\eta_y}\Expc{\sqn{u_y^{t+1}-u_y^\star}_{\boldsymbol{\mathcal{U}}_y}}\\
&\leq \frac{1}{\gamma}\sqn{\mathbf{w}^t-\mathbf{w}^\star}_{\boldsymbol{\mathcal{X}}}+\frac{\gamma}{p^2\eta}\sqn{\mathbf{u}^t-\mathbf{u}^\star}_{\boldsymbol{\mathcal{U}}}+\frac{\gamma}{p^2\eta_y}\sqn{u_y^t-u_y^\star}_{\boldsymbol{\mathcal{U}}_y}\\
&\quad-\gamma \sqn{D^* (\mathbf{u}^t-\mathbf{u}^\star)+D_y^* (u_y^t-u_y^\star)}_{\boldsymbol{\mathcal{X}}}\\
&\quad+\frac{-2p\rho+p\rho^2(2+ \oma)+p\eta(1+2\omega)}{\gamma}\sqn{D\mathbf{\hat{x}}^t}_{\boldsymbol{\mathcal{U}}}\\
&\quad+\frac{-2p\rho_y+p\rho_y^2(2+\omega_s)+p\rho^2(\omega_s+ \oma)+p\eta(\omega_s+2\omega)+p\eta_y(1+\omega_s)}{\gamma}\sqn{D_y\mathbf{\hat{x}}^t}_{\boldsymbol{\mathcal{U}}_y}.
\end{align*}
We need to choose $\rho,\rho_y,\eta,\eta_y$ small enough to remove the last two squared norm terms. 
We choose 
\begin{align*}
\rho=\rho_y&=\frac{1}{2+\oma+2\omega_s},\\
\eta=\eta_y&=\frac{1}{(1+2\omega+2\omega_s)(2+\oma+2\omega_s)}.
\end{align*}
This way,
\begin{align*}
-2p\rho_y+p\rho_y^2(2+\omega_s)+p\rho^2(\omega_s+ \oma)+p\eta(\omega_s+2\omega)+p\eta_y(1+\omega_s)=0
\end{align*}
and
\begin{align*}
-2p\rho+p\rho^2(2+ \oma)+p\eta(1+2\omega)&\leq -2p\rho+p\rho^2(2+\oma+2\omega_s)+p\eta(1+2\omega+2\omega_s)=0.
\end{align*}
Then, according to \citet[Lemma 1]{con22rp},
\begin{align}
\sqnx{\mathbf{w}^t-\mathbf{w}^\star}&=
\sqnx{(\mathrm{Id}-\gamma\nabla \mathbf{f})\mathbf{x}^t-(\mathrm{Id}-\gamma\nabla \mathbf{f})\mathbf{x}^\star} \notag\\
&\leq \max(1-\gamma\mu,\gamma L-1)^2 \sqnx{\mathbf{x}^t-\mathbf{x}^\star}.\label{eqle1}
\end{align}
Moreover, we define the concatenated operator $D_c : \mathbf{x}\in\boldsymbol{\mathcal{X}}\mapsto (D\mathbf{x},D_y\mathbf{x})\in \boldsymbol{\mathcal{U}}\times \boldsymbol{\mathcal{U}}_y$. $D_c$ has full range, since for every $(x_1,\ldots,x_n,y)\in \boldsymbol{\mathcal{U}}\times \boldsymbol{\mathcal{U}}_y$, 
$D_c(x_1,\ldots,x_n,0,y)=(x_1,\ldots,x_n,y)$. Equivalently, $D^*_c$ is injective. For every $(\mathbf{u},u_y)\in \boldsymbol{\mathcal{U}}\times \boldsymbol{\mathcal{U}}_y$, we have
$\sqn{D^* \mathbf{u}+D_y^* u_y}_{\boldsymbol{\mathcal{X}}}=\sqn{D^*_c( \mathbf{u},u_y)}_{\boldsymbol{\mathcal{X}}}\geq \lambda_{\min}(D_cD^*_c)\big(
\sqn{\mathbf{u}}_{\boldsymbol{\mathcal{U}}}+\sqn{u_y}_{\boldsymbol{\mathcal{U}}_y}\big)$, where $\lambda_{\min}(D_cD^*_c)$ is the smallest eigenvalue of $D_cD^*_c$, which is positive because $D^*_c$ is injective. This eigenvalue is 1, with corresponding eigenvectors of the form $(u_1,\ldots,u_n,-\frac{1}{n}\sum_{i=1}^n u_i)$, that form a space of dimension $n$ (the $(n+1)$-th eigenvalue is 2, as mentioned above). Therefore, 
\begin{align*}
\sqn{D^* (\mathbf{u}^t-\mathbf{u}^\star)+D_y^* (u_y^t-u_y^\star)}_{\boldsymbol{\mathcal{X}}}&\geq \sqn{\mathbf{u}^t-\mathbf{u}^\star}_{\boldsymbol{\mathcal{U}}}+\sqn{u_y^t-u_y^\star}_{\boldsymbol{\mathcal{U}}_y}.
\end{align*}
Hence,
\begin{align*}
&\frac{1}{\gamma}\Expc{\sqn{\mathbf{x}^{t+1}-\mathbf{x}^\star}_{\boldsymbol{\mathcal{X}}}}+\frac{\gamma}{p^2\eta}\Expc{\sqn{\mathbf{u}^{t+1}-\mathbf{u}^\star}_{\boldsymbol{\mathcal{U}}}}+
\frac{\gamma}{p^2\eta_y}\Expc{\sqn{u_y^{t+1}-u_y^\star}_{\boldsymbol{\mathcal{U}}_y}}\\
&\leq \frac{1}{\gamma}\max(1-\gamma\mu,\gamma L-1)^2 \sqnx{\mathbf{x}^t-\mathbf{x}^\star}
+\frac{\gamma}{p^2\eta}\sqn{\mathbf{u}^t-\mathbf{u}^\star}_{\boldsymbol{\mathcal{U}}}+\frac{\gamma}{p^2\eta_y}\sqn{u_y^t-u_y^\star}_{\boldsymbol{\mathcal{U}}_y}\\
&\quad-\gamma\sqn{\mathbf{u}^t-\mathbf{u}^\star}_{\boldsymbol{\mathcal{U}}}-\gamma\sqn{u_y^t-u_y^\star}_{\boldsymbol{\mathcal{U}}_y}\\
&= \frac{1}{\gamma}\max(1-\gamma\mu,\gamma L-1)^2 \sqnx{\mathbf{x}^t-\mathbf{x}^\star}
+\frac{\gamma}{p^2\eta}\big(1-p^2\eta\big)\sqn{\mathbf{u}^t-\mathbf{u}^\star}_{\boldsymbol{\mathcal{U}}}\\
&\quad+\frac{\gamma}{p^2\eta_y}\big(1-p^2\eta_y\big)\sqn{u_y^t-u_y^\star}_{\boldsymbol{\mathcal{U}}_y},
\end{align*}
so that
\begin{align}
\Exp{\Psi^{t+1}\;|\;\mathcal{F}^t}&\leq \max\left((1-\gamma\mu)^2,(1-\gamma L)^2,1-p^2\eta\right)\Psi^t.\label{eqrec2b}
\end{align}
Using the tower rule, we can unroll the recursion in \eqref{eqrec2b} to obtain the unconditional expectation of $\Psi^{t+1}$. 

Using classical results on supermartingale convergence \citep[Proposition A.4.5]{ber15}, it follows from \eqref{eqrec2b} that $\Psi^t\rightarrow 0$ almost surely. Almost sure convergence of $\mathbf{x}^t$, $ \mathbf{u}^t$, $u_y^t$ follows. \medskip

Finally, let us consider the general case $k\in [d]$. We observe that the analysis above is separable with respect to the $d$ coordinates of the vectors. From the perspective of a given coordinate $j\in[d]$, either the vector values at this coordinate are updated using communicated information, which happens if $\theta^t=1$ and $j\in\Omega^t$, or they are updated using local information only, which happens if $\theta^t=0$ or $j\notin\Omega^t$. So, since $j\in\Omega^t$ happens with probability $k/d$, the whole analysis above applies with $p$ replaced by $\frac{pk}{d}$.

\section{Proof of Theorem~\ref{theogc}}

 For every $t\geq 0$, we define  $\mathcal{D}_{f_i}^t\eqdef  \mathcal{D}_{f_i}(x_i^t,x^\star)$, for every $i\in[n]$, $\mathcal{D}_{f_s}^t\eqdef  \mathcal{D}_{f_s}(x_s^t,x^\star)$, $\mathcal{D}_{g}^t\eqdef  \mathcal{D}_{g}(y^t,x^\star)$.
 
We follow the same derivations as in Appendix~\ref{secalg2}. However, instead of \eqref{eqle1}, we use
\begin{align}
\sqnx{\mathbf{w}^t-\mathbf{w}^\star}&=
\sqnx{\mathbf{x}^t-\mathbf{x}^\star}-2\gamma \left\langle \nabla \mathbf{f}(\mathbf{x}^t)-\nabla \mathbf{f}(\mathbf{x}^\star),\mathbf{x}^t-\mathbf{x}^\star\right\rangle_{\boldsymbol{\mathcal{X}}}+\gamma^2\sqnx{\nabla \mathbf{f}(\mathbf{x}^t)-\nabla \mathbf{f}(\mathbf{x}^\star)}.
\end{align}
Also, we choose 
\begin{align*}
\rho=\rho_y&=\frac{1}{2+\oma+2\omega_s},\\
\eta=\eta_y&<\frac{1}{(1+2\omega+2\omega_s)(2+\oma+2\omega_s)}
\end{align*}
(for instance $\eta=\eta_y=\frac{0.99}{(1+2\omega+2\omega_s)(2+\oma+2\omega_s)}$). This way,
\begin{align*}
\chi_y\eqdef 2\rho_y-\rho_y^2(2+\omega_s)-\rho^2(\omega_s+ \oma)-\eta(\omega_s+2\omega)-\eta_y(1+\omega_s)>0
\end{align*}
and
\begin{align*}
\chi\eqdef 2\rho-\rho^2(2+ \oma)-\eta(1+2\omega)&\geq 2\rho-\rho^2(2+\oma+2\omega_s)-\eta(1+2\omega+2\omega_s)>0.
\end{align*}
Hence,
\begin{align*}
&\frac{1}{\gamma}\Expc{\sqn{\mathbf{x}^{t+1}-\mathbf{x}^\star}_{\boldsymbol{\mathcal{X}}}}+\frac{\gamma}{p_{t+1}^2\eta}\Expc{\sqn{\mathbf{u}^{t+1}-\mathbf{u}^\star}_{\boldsymbol{\mathcal{U}}}}+
\frac{\gamma}{p_{t+1}^2\eta}\Expc{\sqn{u_y^{t+1}-u_y^\star}_{\boldsymbol{\mathcal{U}}_y}}\\
&\leq \frac{1}{\gamma}\sqnx{\mathbf{x}^t-\mathbf{x}^\star}-2 \left\langle \nabla \mathbf{f}(\mathbf{x}^t)-\nabla \mathbf{f}(\mathbf{x}^\star),\mathbf{x}^t-\mathbf{x}^\star\right\rangle_{\boldsymbol{\mathcal{X}}}+\gamma\sqnx{\nabla \mathbf{f}(\mathbf{x}^t)-\nabla \mathbf{f}(\mathbf{x}^\star)}\\
&\quad+\frac{\gamma}{p_{t+1}^2\eta}\big(1-p_{t+1}^2\eta\big)\sqn{\mathbf{u}^t-\mathbf{u}^\star}_{\boldsymbol{\mathcal{U}}}+\frac{\gamma}{p_{t+1}^2\eta}\big(1-p_{t+1}^2\eta\big)\sqn{u_y^t-u_y^\star}_{\boldsymbol{\mathcal{U}}_y}\\
&\quad-\frac{p_{t+1}\chi}{\gamma}\sqn{D\mathbf{\hat{x}}^t}_{\boldsymbol{\mathcal{U}}}-\frac{p_{t+1}\chi_y}{\gamma}\sqn{D_y\mathbf{\hat{x}}^t}_{\boldsymbol{\mathcal{U}}_y}.
\end{align*}
So, assuming for simplicity that $\gamma\leq \frac{1}{L}$, we have
\begin{align*}
&\frac{1}{\gamma}\Expc{\sqn{\mathbf{x}^{t+1}-\mathbf{x}^\star}_{\boldsymbol{\mathcal{X}}}}+\frac{\gamma}{p_{t+1}^2\eta}\Expc{\sqn{\mathbf{u}^{t+1}-\mathbf{u}^\star}_{\boldsymbol{\mathcal{U}}}}+
\frac{\gamma}{p_{t+1}^2\eta}\Expc{\sqn{u_y^{t+1}-u_y^\star}_{\boldsymbol{\mathcal{U}}_y}}\\
&\leq \frac{1}{\gamma}\sqnx{\mathbf{x}^t-\mathbf{x}^\star}-2\sum_{i=1}^n \mathcal{D}_{f_i}^t - 4n \mathcal{D}_{f_s}^t -2n \mathcal{D}_{g}^t
+\left(\gamma-\frac{1}{L}\right)\sqnx{\nabla \mathbf{f}(\mathbf{x}^t)-\nabla \mathbf{f}(\mathbf{x}^\star)}\\
&\quad-\frac{p_{t+1}\chi}{\gamma}\sqn{D\mathbf{\hat{x}}^t}_{\boldsymbol{\mathcal{U}}}-\frac{p_{t+1}\chi_y}{\gamma}\sqn{D_y\mathbf{\hat{x}}^t}_{\boldsymbol{\mathcal{U}}_y}\\
&\quad+\frac{\gamma}{p_{t+1}^2\eta}\big(1-p_{t+1}^2\eta\big)\sqn{\mathbf{u}^t-\mathbf{u}^\star}_{\boldsymbol{\mathcal{U}}}+\frac{\gamma}{p_{t+1}^2\eta}\big(1-p_{t+1}^2\eta\big)\sqn{u_y^t-u_y^\star}_{\boldsymbol{\mathcal{U}}_y}\\
&\leq \frac{1}{\gamma}\sqnx{\mathbf{x}^t-\mathbf{x}^\star}-2\sum_{i=1}^n \mathcal{D}_{f_i}^t - 4n \mathcal{D}_{f_s}^t -2n \mathcal{D}_{g}^t
-\frac{p_{t+1}\chi}{\gamma}\sqn{D\mathbf{\hat{x}}^t}_{\boldsymbol{\mathcal{U}}}-\frac{p_{t+1}\chi_y}{\gamma}\sqn{D_y\mathbf{\hat{x}}^t}_{\boldsymbol{\mathcal{U}}_y}\\
&\quad+\gamma\left(\frac{1}{p_{t+1}^2\eta}-1\right)\sqn{\mathbf{u}^t-\mathbf{u}^\star}_{\boldsymbol{\mathcal{U}}}+\gamma\left(\frac{1}{p_{t+1}^2\eta}-1\right)\sqn{u_y^t-u_y^\star}_{\boldsymbol{\mathcal{U}}_y}.
\end{align*}
We choose 
\begin{equation*}
p_t = \sqrt{\frac{b}{a+t}}
\end{equation*}
for some  $b\geq \frac{1}{\eta}$ and $a\geq b-1$ (so that $p_t \in (0,1]$ for every $t\geq 1$). Then we have
\begin{equation*}
\frac{1}{p^2_{t+1}\eta}-1 = \frac{a-b\eta+t+1}{b\eta}\leq \frac{a+t}{b\eta}=\frac{1}{p^2_{t}\eta}.
\end{equation*}
Hence, 
\begin{align*}
&\frac{1}{\gamma}\Expc{\sqn{\mathbf{x}^{t+1}-\mathbf{x}^\star}_{\boldsymbol{\mathcal{X}}}}+\frac{\gamma}{p_{t+1}^2\eta}\Expc{\sqn{\mathbf{u}^{t+1}-\mathbf{u}^\star}_{\boldsymbol{\mathcal{U}}}}+
\frac{\gamma}{p_{t+1}^2\eta}\Expc{\sqn{u_y^{t+1}-u_y^\star}_{\boldsymbol{\mathcal{U}}_y}}\\
&\leq \frac{1}{\gamma}\sqnx{\mathbf{x}^t-\mathbf{x}^\star}-2\sum_{i=1}^n \mathcal{D}_{f_i}^t - 4n \mathcal{D}_{f_s}^t -2n \mathcal{D}_{g}^t
-\frac{p_{t+1}\chi}{\gamma}\sqn{D\mathbf{\hat{x}}^t}_{\boldsymbol{\mathcal{U}}}-\frac{p_{t+1}\chi_y}{\gamma}\sqn{D_y\mathbf{\hat{x}}^t}_{\boldsymbol{\mathcal{U}}_y}
\\
&\quad+\frac{\gamma}{p_{t}^2\eta}\sqn{\mathbf{u}^t-\mathbf{u}^\star}_{\boldsymbol{\mathcal{U}}}+\frac{\gamma}{p_{t}^2\eta}\sqn{u_y^t-u_y^\star}_{\boldsymbol{\mathcal{U}}_y}.
\end{align*}
By unrolling the recursion, we have, for every $t\geq 1$,
\begin{align*}
&\frac{1}{\gamma}\Expc{\sqn{\mathbf{x}^{t}-\mathbf{x}^\star}_{\boldsymbol{\mathcal{X}}}}+\frac{\gamma}{p_{t}^2\eta}\Expc{\sqn{\mathbf{u}^{t}-\mathbf{u}^\star}_{\boldsymbol{\mathcal{U}}}}+
\frac{\gamma}{p_{t}^2\eta}\Expc{\sqn{u_y^{t}-u_y^\star}_{\boldsymbol{\mathcal{U}}_y}}\\
&\leq \frac{1}{\gamma}\sqnx{\mathbf{x}^0-\mathbf{x}^\star}+\frac{\gamma a}{\eta b}\sqn{\mathbf{u}^0-\mathbf{u}^\star}_{\boldsymbol{\mathcal{U}}}+\frac{\gamma a}{ \eta b}\sqn{u_y^0-u_y^\star}_{\boldsymbol{\mathcal{U}}_y}\\
&\quad- \sum_{t'=0}^{t-1}\Exp{2\sum_{i=1}^n \mathcal{D}_{f_i}^{t'} +4n \mathcal{D}_{f_s}^{t'} +2n \mathcal{D}_{g}^{t'}
+\frac{p_{t'+1}\chi}{\gamma}\sqn{D\mathbf{\hat{x}}^{t'}}_{\boldsymbol{\mathcal{U}}}+\frac{p_{t'+1}\chi_y}{\gamma}\sqn{D_y\mathbf{\hat{x}}^{t'}}_{\boldsymbol{\mathcal{U}}_y}
}.
\end{align*}
This implies that 
\begin{align*}
&\sum_{t=0}^{\infty}\Exp{2\sum_{i=1}^n \mathcal{D}_{f_i}^{t} +4n \mathcal{D}_{f_s}^{t} +2n \mathcal{D}_{g}^{t}
+\frac{p_{t+1}\chi}{\gamma}\sqn{D\mathbf{\hat{x}}^{t}}_{\boldsymbol{\mathcal{U}}}+\frac{p_{t+1}\chi_y}{\gamma}\sqn{D_y\mathbf{\hat{x}}^{t}}_{\boldsymbol{\mathcal{U}}_y}}\\
&\leq \frac{1}{\gamma}\sqnx{\mathbf{x}^0-\mathbf{x}^\star}+\frac{\gamma a}{\eta b}\sqn{\mathbf{u}^0-\mathbf{u}^\star}_{\boldsymbol{\mathcal{U}}}+\frac{\gamma a}{ \eta b}\sqn{u_y^0-u_y^\star}_{\boldsymbol{\mathcal{U}}_y},
\end{align*}
so that $\Exp{2\sum_{i=1}^n \mathcal{D}_{f_i}^{t} +4n \mathcal{D}_{f_s}^{t} +2n \mathcal{D}_{g}^{t}
+\frac{p_{t+1}\chi}{\gamma}\sqn{D\mathbf{\hat{x}}^{t}}_{\boldsymbol{\mathcal{U}}}+\frac{p_{t+1}\chi_y}{\gamma}\sqn{D_y\mathbf{\hat{x}}^{t}}_{\boldsymbol{\mathcal{U}}_y}}\rightarrow 0$ as $t\rightarrow +\infty$. 
Moreover, for every $t\geq 0$,
\begin{align*}
&\Expc{\sqn{\mathbf{u}^{t}-\mathbf{u}^\star}_{\boldsymbol{\mathcal{U}}}}+
\Expc{\sqn{u_y^{t}-u_y^\star}_{\boldsymbol{\mathcal{U}}_y}}\\
&\leq p_t^2\left(\frac{\eta}{\gamma^2}\sqnx{\mathbf{x}^0-\mathbf{x}^\star}+\frac{ a}{ b}\sqn{\mathbf{u}^0-\mathbf{u}^\star}_{\boldsymbol{\mathcal{U}}}+\frac{ a}{  b}\sqn{u_y^0-u_y^\star}_{\boldsymbol{\mathcal{U}}_y}\right)\\
& \frac{b}{a+t}\left(\frac{\eta}{\gamma^2}\sqnx{\mathbf{x}^0-\mathbf{x}^\star}+\frac{ a}{ b}\sqn{\mathbf{u}^0-\mathbf{u}^\star}_{\boldsymbol{\mathcal{U}}}+\frac{ a}{  b}\sqn{u_y^0-u_y^\star}_{\boldsymbol{\mathcal{U}}_y}\right)\\
&=\mathcal{O}\left(\frac{1}{t}\right).
\end{align*}
Furthermore, for every 
$T\geq 1$, let $\tilde{t}$ be chosen uniformly at random in $\{0,\ldots,T-1\}$. Then
\begin{align*}
&\Exp{2\sum_{i=1}^n \mathcal{D}_{f_i}^{\tilde{t}} +4n \mathcal{D}_{f_s}^{\tilde{t}} +2n \mathcal{D}_{g}^{\tilde{t}}
+\frac{p_{\tilde{t}+1}\chi}{\gamma}\sqn{D\mathbf{\hat{x}}^{\tilde{t}}}_{\boldsymbol{\mathcal{U}}}+\frac{p_{\tilde{t}+1}\chi_y}{\gamma}\sqn{D_y\mathbf{\hat{x}}^{\tilde{t}}}_{\boldsymbol{\mathcal{U}}_y}}\\
&=\frac{1}{T}\sum_{t=0}^{T-1} \Exp{2\sum_{i=1}^n \mathcal{D}_{f_i}^{t} +4n \mathcal{D}_{f_s}^{t} +2n \mathcal{D}_{g}^{t}
+\frac{p_{t+1}\chi}{\gamma}\sqn{D\mathbf{\hat{x}}^{t}}_{\boldsymbol{\mathcal{U}}}+\frac{p_{t+1}\chi_y}{\gamma}\sqn{D_y\mathbf{\hat{x}}^{t}}_{\boldsymbol{\mathcal{U}}_y}}\\
&\leq \frac{1}{T}\left( \frac{1}{\gamma}\sqnx{\mathbf{x}^0-\mathbf{x}^\star}+\frac{\gamma a}{\eta b}\sqn{\mathbf{u}^0-\mathbf{u}^\star}_{\boldsymbol{\mathcal{U}}}+\frac{\gamma a}{ \eta b}\sqn{u_y^0-u_y^\star}_{\boldsymbol{\mathcal{U}}_y}\right).
\end{align*}
Thus, given $\epsilon>0$, by choosing 
\begin{align*}
T&\geq \frac{1}{2\epsilon}\left( \frac{1}{\gamma}\sqnx{\mathbf{x}^0-\mathbf{x}^\star}+\frac{\gamma a}{\eta b}\sqn{\mathbf{u}^0-\mathbf{u}^\star}_{\boldsymbol{\mathcal{U}}}+\frac{\gamma a}{ \eta b}\sqn{u_y^0-u_y^\star}_{\boldsymbol{\mathcal{U}}_y}\right),
\end{align*}
we have 
\begin{align*}
&\Exp{2\sum_{i=1}^n \mathcal{D}_{f_i}^{\tilde{t}} +4n \mathcal{D}_{f_s}^{\tilde{t}} +2n \mathcal{D}_{g}^{\tilde{t}}
+\frac{p_{\tilde{t}+1}\chi}{\gamma}\sqn{D\mathbf{\hat{x}}^{\tilde{t}}}_{\boldsymbol{\mathcal{U}}}+\frac{p_{\tilde{t}+1}\chi_y}{\gamma}\sqn{D_y\mathbf{\hat{x}}^{\tilde{t}}}_{\boldsymbol{\mathcal{U}}_y}}\leq 2\epsilon.
\end{align*}
In particular,
\begin{align}
&\Exp{\sum_{i=1}^n \mathcal{D}_{f_i}^{\tilde{t}} +2n \mathcal{D}_{f_s}^{\tilde{t}} +n \mathcal{D}_{g}^{\tilde{t}}
}\leq \epsilon.\label{eqe1}
\end{align}
Also, since $(p_t)_t$ is decreasing,
\begin{align*}
\frac{p_{T}\chi}{\gamma}\sqn{D\mathbf{\hat{x}}^{\tilde{t}}}_{\boldsymbol{\mathcal{U}}}+\frac{p_{T}\chi_y}{\gamma}\sqn{D_y\mathbf{\hat{x}}^{\tilde{t}}}_{\boldsymbol{\mathcal{U}}_y}&\leq \frac{p_{\tilde{t}+1}\chi}{\gamma}\sqn{D\mathbf{\hat{x}}^{\tilde{t}}}_{\boldsymbol{\mathcal{U}}}+\frac{p_{\tilde{t}+1}\chi_y}{\gamma}\sqn{D_y\mathbf{\hat{x}}^{\tilde{t}}}_{\boldsymbol{\mathcal{U}}_y},
\end{align*}
so that 
\begin{align*}
\Exp{\frac{\chi}{\gamma}\sqn{D\mathbf{\hat{x}}^{\tilde{t}}}_{\boldsymbol{\mathcal{U}}}+\frac{\chi_y}{\gamma}\sqn{D_y\mathbf{\hat{x}}^{\tilde{t}}}_{\boldsymbol{\mathcal{U}}_y}}&\leq \frac{2}{\epsilon p_T}=2\epsilon \sqrt{\frac{a+T}{b}},
\end{align*}
and if $T=\Theta\left( \frac{1}{2\epsilon}\left( \frac{1}{\gamma}\sqnx{\mathbf{x}^0-\mathbf{x}^\star}+\frac{\gamma a}{\eta b}\sqn{\mathbf{u}^0-\mathbf{u}^\star}_{\boldsymbol{\mathcal{U}}}+\frac{\gamma a}{ \eta b}\sqn{u_y^0-u_y^\star}_{\boldsymbol{\mathcal{U}}_y}\right)\right)$,
\begin{align*}
\Exp{\frac{\chi}{\gamma}\sqn{D\mathbf{\hat{x}}^{\tilde{t}}}_{\boldsymbol{\mathcal{U}}}+\frac{\chi_y}{\gamma}\sqn{D_y\mathbf{\hat{x}}^{\tilde{t}}}_{\boldsymbol{\mathcal{U}}_y}}&=\mathcal{O}\left(\sqrt{\epsilon}\right).
\end{align*}
Finally, the expectation of the number of communication rounds over the first $T\geq 1$ iterations is 
\begin{equation*}
\sum_{t=1}^{T} p_t = \Theta(\sqrt{T}),
\end{equation*}
so that \eqref{eqe1} is achieved with $\Theta(\sqrt{T})=\Theta\left(\frac{1}{\sqrt{\epsilon}}\right)$ communication rounds.

\end{document}